\documentclass[aop,secthm,seceqn,dvips]{arximspdf}

\doi{10.1214/009117907000000259}
\volume{36}
\issue{1}
\pubyear{2008}
\firstpage{363}
\lastpage{396}

\makeatletter

\newproclaim{remarks}{Remarks}
\newproclaim{remark}{Remark}
\newtheorem{lemm}{Lemma}[section]

\def\descriptionlabel#1{\hspace\labelsep\it #1\ }

\makeatother

\begin{document}
\begin{frontmatter}

\title{Variance asymptotics and central limit theorems for
generalized growth processes with applications to convex hulls\\ and
maximal points}
\runtitle{Variance asymptotics and central limit theorems}
\pdftitle{Variance asymptotics and central limit theorems for
generalized growth processes with applications to convex hulls and
maximal points}

\begin{aug}

\author[A]{\fnms{T.} \snm{Schreiber}\ead[label=e1]{tomeks@mat.uni.torun.pl}\thanksref{t1}}
and
\author[B]{\fnms{J. E.} \snm{Yukich}\corref{}\ead[label=e2]{joseph.yukich@lehigh.edu}\thanksref{t2}}

\thankstext{t1}{Supported in part by Polish Minister
of Scientific Research and Information Technology Grant 1 P03A 018 28
(2005--2007).}
\thankstext{t2}{Supported by NSF Grant DMS-02-03720.}

\runauthor{T. Schreiber and J. E. Yukich}

\affiliation{Nicholas Copernicus University and Lehigh University}

\address[A]{Faculty of Mathematics and Computer Science\\
Nicholas Copernicus University\\
Toru\'n\\
Poland\\
\printead{e1}} 

\address[B]{Department of Mathematics\\
Lehigh University\\
Bethlehem, Pennsylvania 18015\\
USA\\
\printead{e2}}

\end{aug}

\received{\smonth{5} \syear{2006}}
\revised{\smonth{2} \syear{2007}}

%
\begin{abstract}
We show that the random point measures induced by vertices in the
convex hull of a Poisson sample on the unit
ball, when properly scaled and centered, converge to those of a
mean zero Gaussian field. We establish limiting variance and
covariance asymptotics in terms of the density of the Poisson
sample. Similar results hold for the point measures induced by the
maximal points in a Poisson sample. The approach involves
introducing a generalized spatial birth growth process allowing
for cell overlap.
\end{abstract}

%

\end{frontmatter}

\section{Introduction, main results}\label{GenRes}

Given $X_i, i \geq1$, i.i.d. random variables with
values in a $d$-dimensional convex set $S$, $d \geq2,$
a classic problem in convex
geometry involves determining the distribution of the number of
points in the set of extreme points $\mathcal{V} (\{X_i\}_{i=1}^n)$,
defined as the vertices in the convex hull of
$\{X_i\}_{i=1}^n$. This problem was first considered by R\'enyi
and Sulanke \cite{RS}, with recent notable progress by Reitzner
\cite{Re1,Re2,Re3,Re4} and Vu \cite{VV}.

A closely related
problem involves determining, for a given $K \subset\mathbb{R}^d$, the
distribution of the number of points in
the set $\mathcal{M}_K(\{X_i\}_{i=1}^n)$ of $K$-maximal points,
where a point $X_j$ belongs to $\mathcal{M}_K(\{X_i\}_{i=1}^n)$ iff
$(X_j \oplus K) \ \cap\{X_i\}_{i=1}^n = X_j$, 
where here and henceforth, for all $B \subset\mathbb{R}^d$ and $x \in
\mathbb{R}^d$ we write $x \oplus B := \{x + y\dvtx y \in B\}$. When
$K$ is
$(\mathbb{R}_+)^d$, then $\mathcal{M}_K(\{X_i\}_{i=1}^n)$ is simply the
set of
maximal points, that is, those points $X_j$ in $\{X_i\}_{i=1}^n$
having the property that no point $X_i$, $i \neq j$, exceeds it
in all coordinates. The limit theory for the number of maximal
points in $\mathcal{M}_K(\{X_i\}_{i=1}^n)$ was first considered by
R\'enyi \cite{Re} and Barndorff-Nielsen and Sobel \cite{BS}. Chen,
Hwang and Tsai \cite{Ch} surveys the vast literature, which includes books
by Ehrgott \cite{Eh}, Pomerol and Barba-Romero \cite{PB}, and
recent papers of \cite{BCHL,BHLT,BX,BY4,De}.

In this paper we establish convergence of the finite-dimensional
distributions of the re-scaled point measures induced by the
random point sets $\mathcal{V} ({\mathcal{P}_{{\lambda}\rho} })$, where
${\mathcal{P}_{{\lambda}\rho} }$ denotes a Poisson
point process of intensity ${\lambda}{\rho}$ on $B_d$, the unit radius
\mbox{$d$-dimensional} ball centered at the origin
and where ${\rho}$ is a continuous density on $B_d$. For
sets $K := \{(w_1,\ldots,w_d) \dvtx w_d \geq(w_1^2 + \cdots+
w_{d-1}^2)^{\alpha/2}\}$, where $\alpha\in(0,1]$ is fixed, we
also establish convergence of the finite-dimensional distributions
of the point measures induced by $\mathcal{M}_K ({\mathcal
{P}_{{\lambda}\rho}
})$, where ${\mathcal{P}_{{\lambda}\rho} }$
denotes the Poisson point process of intensity ${\lambda}{\rho}$ on $A
\times\mathbb{R}_+$, where $A \subset\mathbb{R}^{d-1}$ is compact
and convex and
where ${\rho}\dvtx A \times\mathbb{R}_+$ is continuous.
These results are facilitated by
introducing a generalized spatial birth--growth process as a means
toward obtaining explicit variance asymptotics and central limit
theorems for random measures arising in convex geometry.
The relevant spatial birth--growth process, possibly of independent
interest, modifies
the classical spatial birth--growth process introduced by
Kolmogorov \cite{Ko} as a model for crystal growth by allowing
the possibility of \textit{cell overlap}. As in \cite{Ko}, cells may
grow at nonconstant growth rates.


In the context of the set of extreme points $\mathcal{V} ({\mathcal
{P}_{{\lambda}
\rho} })$, the
approach taken here adds to the work of Reitzner
\cite{Re1,Re2,Re3,Re4} and Vu \cite{VV} in the following ways.
First, the present set-up establishes convergence of the
finite-dimensional distributions of the canonical point measures induced
by $\mathcal{V} ({\mathcal{P}_{{\lambda}\rho} })$, whereas \cite
{Re1,Re2,Re3,Re4} and \cite{VV}
deal with one-dimensional central limit theorems. Second, we
establish a formula for variance and covariance asymptotics.
Third, the present paper concerns the limit theory for nonuniform
samples, whereas \cite{Re1,Re2,Re3,Re4} and \cite{VV} treat
uniform random samples.

In the context of the set of maximal points $\mathcal{M}_K({\mathcal
{P}_{{\lambda}
\rho} })$, the
present set-up establishes convergence of the finite-dimensional
distributions of the canonical point measures induced by $\mathcal
{M}_K({\mathcal{P}_{{\lambda}\rho} })$, with covariances, whereas
previous work
\cite{BX,De}
is concerned with one dimensional central limit theorems without a
formula for covariance asymptotics and/or is limited to the case
when $K$ is a cone \cite{BY4}.

\subsection{Terminology\textup{,} $\psi$-growth processes}\label
{Terminology}

Let the function $\psi\dvtx\mathbb{R}_+ \to\mathbb{R}_+$ satisfy
the following
conditions:
\begin{longlist}
\item[($\Psi$1)] $\psi$ is monotone and $\lim_{l\to\infty}
\psi(l) =
\infty$, and
\item[($\Psi$2)] there exists $\alpha> 0$ such that
$\psi(l) = l^{\alpha} (1+o(1))$ for $l$ small enough.
\end{longlist}
Let $\mathbf{0}$ denote the origin of $\mathbb{R}^{d-1}$, $d
\geq2$, and let
$|y|$ denote the Euclidean norm of $y \in\mathbb{R}^d$.
We define $K[\mathbf{0}]$ to be the $\psi$-epigraph $\{ (y,h) \in
\mathbb{R}^{d-1}
\times\mathbb{R}_+ \dvtx
\ h \geq\psi(|y|) \}$ and, more generally, for $\bar{x} := (x,h_x)
\in\mathbb{R}^{d-1}
\times\mathbb{R}_+$,
we define its \mbox{$\psi$-epigraph} (or upward cone) by
%
%
\begin{equation}\label{upcone}
K[\bar{x}]:= \bar{x} \oplus K[\mathbf{0}]:= \{ (y,h) \in\mathbb{R}^{d-1}
\times\mathbb{R}_+ \dvtx h \geq h_x + \psi(|y-x|) \}.
\end{equation}
%
Given a point set $\mathcal{X}
\subseteq\mathbb{R}^{d-1} \times\mathbb{R}_+,$ a point $\bar{x}
\in\mathcal{X}$ is called
\textit{$\psi$-extremal} in $\mathcal{X}$ iff $K[\bar{x}] \not
\subseteq
\bigcup_{\bar{y}
\in\mathcal{X} \setminus\{ \bar{x} \}} K[\bar{y}]$, that is to
say the
$\psi$-epigraph
of $\bar{x}$ is not completely covered by the union of the $\psi
$-epigraphs of points in
$\mathcal{X} \setminus\{x\}$.
Define the functional
%
%
\begin{equation}\label{XiPsi}
\xi(\bar{x},\mathcal{X}):=
\xi(\psi; \bar{x},\mathcal{X}):= \cases{
1,&\quad if $\bar{x}$ is $\psi$-extremal in
$\mathcal{X}$, \cr
0,&\quad otherwise.
}
\end{equation}
With $D$ standing for some bounded domain
in $\mathbb{R}^{d-1} \times\mathbb{R}_+,$ we consider the version
$\xi_D(\cdot,\cdot)$ of
$\xi(\cdot,\cdot)$ restricted to $D,$ by setting $\xi_D(\bar
{x},\mathcal{X})$ to be $1$ iff $K[\bar{x}] \cap D \not\subseteq
\bigcup_{\bar{y} \in
(\mathcal{X}\setminus\{ \bar{x}\}) \cap D} K[\bar{y}],$ in which
case we declare $\bar{x}$ to be $\psi$-extremal in $D \cap\mathcal{X},$
and otherwise we set $\xi_D(\bar{x},\mathcal{X})$ to be zero.
In case $\bar{x} \notin\mathcal{X}$ we abbreviate notation and write
$\xi(\bar{x},\mathcal{X})$ for $\xi(\bar{x},\mathcal{X} \cup\bar{x})$
and similarly for $\xi_D(\bar{x},\mathcal{X})$.

To provide a physical interpretation
of these functionals, we regard $\mathbb{R}^{d-1} \times\mathbb{R}_+$
as $d$-dimensional
space time, with $\mathbb{R}_+$ standing for the time coordinate,
and we interpret the graph $\partial(K[\bar{x}]),\; \bar{x}:= (x,t)$,
as the
boundary of a $(d-1)$-dimensional spherical particle born at
$x$ at time $t$
(at which time it has initial radius zero) and growing thereupon with
radial speed
$v(t):= \frac{d}{dt}[\psi^{-1}(t)]$, provided the derivative
exists. The particles (spheres) grow independently and
do not exhibit exclusion, that is, they
may overlap or penetrate one another.
A particle is \textit{extreme} iff at some time it is not
completely covered by other particles.
When $\psi$ is the identity, so that the $\psi$ graph
gives a cone, we see that $\psi$-extremal points coincide with
maximal points \cite{BY4}.

In the context of this representation,
it should be noted that, unlike the one stated here, the classic growth
process (see, e.g., \cite{BY2,CQ,Ko,PY2}) assumes that
particles, upon being born at random locations $x \in\mathbb{R}^{d-1}$
at random times $h_x \in\mathbb{R}^+$, form a cell by growing radially
in all directions with a possibly nonconstant speed, that is, with
$\psi$ possibly nonlinear. When one growing cell touches
another, it stops growing in that direction,
that is, no overlap is allowed. Furthermore, a particle born
inside an existing cell is \textit{discarded}, otherwise it is
\textit{accepted}. Letting $\hat{\xi}(\bar{x},\mathcal{X})$ be
zero or one
according to whether $\bar{x}$ is accepted or not, this paper
also considers such functionals $\hat{\xi}$.


The growth process giving rise to the functional $\xi$ will henceforth be
called the \textit{$\psi$-growth process with overlap},
while the process corresponding to
$\hat{\xi}$ will be referred to as the \textit{$\psi$-growth process
without overlap}.
This paper will mainly concentrate on applications
of the first concept and the corresponding functional $\xi,$ but
the subsequently developed general theory also treats the latter
concept in the special case of linear $\psi$.
Throughout, let $A$ be a compact convex subset of $\mathbb{R}^{d-1}.$
We shall also admit the case $A := \mathbb{R}^{d-1}$ in the sequel, in which
case we assume that $\rho$ is uniformly bounded.
Consider a density function $\rho$ on $A_+ := A \times\mathbb{R}_+,$
not necessarily
integrable, such that
\begin{longlist}
\item[(R1)] $\rho$ is continuous on $A_+,$
\item[(R2)] there exists a constant $\delta\geq0$ and a continuous
function $\rho_0\dvtx A \to\mathbb{R}_+$ bounded away from zero such that
\[
\rho(x,h) = \rho_0(x) h^{\delta}
\bigl(1+o(1)\bigr)
\]
for $h$ small enough and
$\rho(x,h) = O (h^{\delta})$ for large $h$
uniformly in $x \in A$.
\end{longlist}



For $\lambda> 0$, we recall that ${\mathcal{P}}_{{\lambda}\rho}$
denotes the Poisson
point process on $A_+$
with intensity measure ${\lambda}\rho(x,h) \,dx\, dh$. The ``extreme
point'' empirical
measures $\mu^{\xi}_{{\lambda}\rho}$ and $\mu^{\hat{\xi
}}_{{\lambda}\rho}$
generated by
${\mathcal{P}}_{{\lambda}\rho}$ are
%
%
\begin{equation}\label{Poiss}
\mu^{\xi}_{{\lambda}\rho} := \sum_{\bar{x} \in
{\mathcal{P}}_{{\lambda}\rho}}
\xi(\bar{x},{\mathcal{P}}_{{\lambda}\rho}) \delta_{\bar{x}}
\end{equation}
and
%
%
\begin{equation}\label{Poisshat}
\mu^{\hat\xi}_{{\lambda}\rho} := \sum_{\bar{x} \in{\mathcal
{P}}_{{\lambda}\rho}}
\hat\xi(\bar{x},{\mathcal{P}}_{{\lambda}\rho}) \delta_{\bar
{x}},
\end{equation}
with $\delta_{x}$ standing for the
unit point mass at $x \in\mathbb{R}^d$. For any random measure
$\sigma$
on $\mathbb{R}^d$, we write $\bar{\sigma}$ for its centered version
$\sigma- \mathbb{E}[\sigma]$, so that, for example,
$ \bar\mu^{\xi}_{{\lambda}\rho} := \mu^{\xi}_{{\lambda}\rho} -
\mathbb{E}[ \mu^{\xi}_{{\lambda}\rho}]$.

Notice that for small $\alpha$ the
upward cones $K[\bar{x}]$ have relatively narrow aperatures,
making it less likely that cones having apexes with a small
temporal coordinate get covered by $\psi$-epigraphs, that is, one
expects more $\psi$-extreme points as $\alpha$ gets smaller. Also,
roughly speaking, for small $\delta$, one expects more points in
${\mathcal{P}}_{{\lambda}\rho}$ with small temporal coordinate and
thus more
$\psi$-extreme points in this case as well. One of the goals of
this paper is to show (see Theorem \ref{LLN}) that the expected
total mass of the extreme point empirical measures
(\ref{Poiss})--(\ref{Poisshat}) is asymptotically proportional to
${\lambda}^{\tau}$, where
%
%
\begin{equation}\label{TAUU}
\tau:= \tau(d, \alpha, \delta) := {d-1 \over d -1 + \alpha(1 +
{\delta}) }.
\end{equation}
More general goals include establishing the variance asymptotics
and the convergence of the finite-dimensional distributions of the
appropriately scaled
measures (\ref{Poiss})--(\ref{Poisshat}) to
Gaussian distributions (see Theorems
\ref{VAR} and \ref{CLT}) and to treat the applications
to extreme and maximal points described at the outset.

\textit{Notation.}
Given $\alpha> 0$, put
%
%
\begin{equation}\label{SCALINGLIMIT1}
\psi^{(\infty)}(l) := l^{\alpha}.
\end{equation}
Recalling the definition of $\xi$, we define the
functional $\xi^{(\infty)}$ by $\xi^{(\infty)}(\cdot, \cdot) :=
\xi(\psi^{(\infty)}; \cdot, \cdot)$ and similarly for
$\hat\xi^{(\infty)}.$ We also let $\mathcal{P}_*$ stand for the Poisson
point process
in $\mathbb{R}^{d-1} \times\mathbb{R}_+$ with intensity measure
$h^{\delta}\, dx\, dh$.

For all $\bar{x}:=(x, h_x)$ and $\bar{y}: = (y, h_y)$, let
\[
m^{(\infty)}(\bar{x}):= \mathbb{E}\bigl[\xi^{(\infty)} (\bar
{x},\mathcal{P}_* )\bigr]
\]
and
\begin{eqnarray*}
c_*^{(\infty)}(\bar{x}, \bar{y})&:=& \mathbb{E}\bigl[ \xi
^{(\infty)}
(\bar
{x},\mathcal{P}_* \cup\bar{y}) \xi^{(\infty)} (\bar{y},\mathcal
{P}_* \cup
\bar{x})\bigr]\\
&&{} - \mathbb{E}\bigl[ \xi^{(\infty)} (\bar{x},\mathcal
{P}_*)\bigr]
\mathbb{E}
[\xi^{(\infty)} (\bar{y},\mathcal{P}_* )]
\end{eqnarray*}
respectively denote the one and two point correlation functions
for the $\psi^{(\infty)}$ growth process with overlap.

For sets $A$ and $B \subset\mathbb{R}^d$, let $d(A,B):=
\inf\{|x-y|\dvtx x \in A,\ y \in B\}$. Let $B_d(y,r)$ denote the
$d$-dimensional
Euclidean ball centered at $y \in\mathbb{R}^d$ with radius $r \in(0,
\infty)$.

Given a subset $B$ of $\mathbb{R}^d$, let $\mathcal{C}_b(B)$ denote the
bounded continuous functions on~$B$. For any signed measure $\mu$
on $A_+$ and $f \in\mathcal{C}_b(A_+)$, let $\langle f, \mu\rangle
:= \int f\, du.$ Unless otherwise specified, $C$ denotes a generic positive
constant whose value may change from line to line.

\subsection{Limit theory for $\Psi$-growth functionals}

For all $f \in\mathcal{C}_b(A_+)$ with $A \subset\mathbb{R}^{d-1}$
compact and
convex, we define the average of the product of $f$ and the one and
two point correlation functions as follows:
%
%
\begin{equation}\label{onept}
I(f):= \int_{A} \int_{0}^{\infty} f(x,0)
m^{(\infty)} (\mathbf{0},h') \rho_0^{\tau}(x) (h')^{\delta}\,
dh'\, dx
\end{equation}
and
\begin{eqnarray}\label{twoptt}
J(f)&:=& \int_A \int_{0}^{\infty}
\int_{\mathbb{R}^{d-1}} \int_0^{\infty} f(x,0)
c_*^{(\infty)}((\mathbf{0},h'), (y',h_y'))\nonumber\\[-8pt]\\[-8pt]
&&\phantom{\int_A \int_{0}^{\infty}
\int_{\mathbb{R}^{d-1}} \int_0^{\infty}}
{}\times\rho_0^{\tau}(x)
(h_y')^{\delta} (h')^{\delta} \,d h_y' \,dy' \,dh' \,dx.\nonumber
\end{eqnarray}
The finiteness of $I(f)$ follows by Lemmas \ref{expbds} and
\ref{L1conv} [see the bound (\ref{IFF})], whereas the finiteness of
$J(f)$ follows from Lemmas \ref{corlimit} and \ref{corbds} [see
the bound (\ref{JFF})] which imply rapid enough decay of
two-point correlation functions.

The following are our main results. We state the results for
$\mu^{\xi}_{{\lambda}\rho}$ and note that analogous results hold for
$\mu^{\hat\xi}_{{\lambda}\rho}$ when $\psi$ is linear. %
The first result specifies first-order behavior, whereas the second
provides second-order
asymptotics.

\begin{thm}\label{LLN}
We have for
all $f \in\mathcal{C}_b(A_+)$
%
%
\begin{equation}\label{LLN1} \lim_{{\lambda}\to\infty} {\lambda
}^{-\tau} \mathbb{E}[\langle f,
\mu^{\xi}_{{\lambda}\rho} \rangle] = I(f).
\end{equation}
%
\end{thm}


\begin{thm}\label{VAR}
We have for all $f \in\mathcal{C}_b(A_+)$
%
%
\begin{equation}\label{varlimit} \lim_{{\lambda}\to\infty}
{\lambda}^{-\tau} \operatorname{Var}[\langle f, \mu^{\xi
}_{{\lambda}\rho}
\rangle] = I(f^2) + J(f^2).
\end{equation}
%
\end{thm}


The next result establishes the convergence of the finite-dimensional
distributions of $({\lambda}^{-\tau/2}
\overline{\mu}^{\xi}_{{\lambda}\rho})$.
\begin{thm} \label{CLT}
The finite-dimensional distributions ${\lambda}^{-\tau/2}(\langle f_1,
\bar{\mu}^{\xi}_{{\lambda}} \rangle,\ldots,\break\langle f_k,
\bar{\mu}^{\xi}_{{\lambda}} \rangle), f_1,\ldots,f_k \in
\mathcal{C}_b(A_+),$
of $({\lambda}^{-\tau/2} \bar{\mu}^{\xi}_{{\lambda}\rho})$
converge as
${\lambda}\to\infty$ to those of a mean zero Gaussian field with
covariance kernel
%
%
\begin{equation}\label{CLT1}
(f,g) \mapsto I(fg) + J(fg), \qquad f, g \in
\mathcal{C}_b(A_+).
\end{equation}
%
\end{thm}

Section \ref{s2} describes applications of $\psi$-growth processes
with overlap, as given by the general limits of Theorems
\ref{LLN}--\ref{CLT}, to convex hulls and maximal points of
i.i.d. samples.

\begin{remarks*}
(i) \textit{Applications to the $\psi$-growth process
$\hat\xi$ without overlap}. The results of Theorems \ref{LLN}--\ref{CLT}
for the functional $\hat\xi$ provide variance
asymptotics and central limit theorems for the classic spatial
birth--growth model in $\mathbb{R}^{d-1}$, whereby seeds are born at random
locations in $ \mathbb{R}^{d-1}$ and times in $\mathbb{R}_+$
according to the
Poisson point process ${\lambda}{\mathcal{P}}_{{\lambda}\rho}$ on
${\lambda}^{1/d}A \times
\mathbb{R}_+$ and grow linearly in time.
Theorems \ref{LLN}--\ref{CLT} for $\hat\xi$ provide a central
limit theorem for the number of seeds accepted in such models.
This generalizes and extends \cite{BY2,PY2}, which builds on
work of Chiu and Quine \cite{CQ,CQa}, Chiu \cite{Chiu} and
Chiu and
Lee \cite{CL}, which do not consider convergence of
finite-dimensional distributions and which often restrict to models with
homogeneous temporal input.
\smallskipamount=0pt
\begin{longlist}[(iii)]
\item[(ii)] \textit{Scaling.} The scaling ${\lambda}^{-\tau}$
arises in the following way. From a conceptual and analytic
point of view, it is convenient to re-scale the $\psi$-growth
process in time and space so as to obtain an equivalent growth
process on Poisson points of approximately unit intensity density
on a region of volume ${\lambda}$. 
The scaling is designed to asymptotically preserve the
$\psi$-epigraphs and the behavior of the density locally close
to $h = 0$.

To achieve
this, we scale $A_+$ in the $d-1$ spatial directions by ${\lambda
}^{{\beta}}$
and in the temporal direction by ${\lambda}^{\gamma}$. Under this
temporal scaling and under (R2), the density $\rho$ exhibits
growth $(h {\lambda}^{\gamma})^{\delta}$ for small temporal $h$, and we
thus require ${\lambda}^{{\beta}(d-1) + \gamma(1 + \delta)} =
{\lambda}$. This
scaling maps $|x|$ and $h_x$ to ${\lambda}^{{\beta}}|x|$ and
${\lambda}^{\gamma}
h_x$, respectively, and therefore, it asymptotically preserves the
$\psi$-epigraphs and condition ($\Psi$2), provided
$({\lambda}^{{\beta}}|x|)^ \alpha= {\lambda}^{\gamma}h_x (1+o(1))$ for
$(x,h_x)$ lying on the graph of $\psi$, that is, $h_x = \psi(x).$
Since $h_x = |x|^{\alpha}
(1+o(1))$ for such $(x,h_x),$
we require ${\lambda}^{{\beta}\alpha} =
{\lambda}^{\gamma}$. We thus require the relations
\[
\beta(d-1) +
\gamma(1 + \delta) = 1 \quad\mbox{and} \quad \beta\alpha=
\gamma,
\]
which yields these values for the scaling exponents
%
%
\begin{equation}
\label{twodef}
\beta= \frac{\gamma} {\alpha} \quad\mbox{and} \quad
\gamma= \frac{\alpha}{(d-1) + \alpha(1 + \delta)}.
\end{equation}
Given the re-scaled $\psi$-growth process on
${\lambda}^{\beta} A \times\mathbb{R}_+$, we expect that a point is
\mbox{$\psi$-extremal} (i.e., $\xi= 1$) iff its time coordinate is
small. Thus, the functional $\mu_{{\lambda}\rho}^{\xi}(A_+)$ should
exhibit growth proportional to the Lebesgue measure of ${\lambda
}^{{\beta}}
A$, that is, proportional to ${\lambda}^{{\beta}(d-1)} = {\lambda
}^{\tau}.$ In the
special case when $\delta= 0$ and the growth is linear ($\alpha=
1$) the $\psi$-epigraphs are preserved by time and space scaling
by ${\lambda}^{1/d}$, that is, $\gamma= 1/d = \beta$. Thus, $\tau=
(d-1)/d$ in this case.



\item[(iii)] \textit{de-Poissonization.} In Section \ref
{ApplSection} we de-Poissonize
Theorems \ref{LLN}--\ref{CLT} when $\alpha\in(0,1]$. In other
words, we obtain the identical limit theory when ${\mathcal
{P}}_{{\lambda
}\rho}$ is
replaced by i.i.d. random variable $X_1,\ldots,X_n$, chosen in
$A_+$ according to
the density $\rho,$ assumed to be integrable to $1$.
We expect similar de-Poissonization results for $\alpha
> 1$, but are unable to prove this.

\item[(iv)] We have not tried to establish a.s. convergence in
(\ref{LLN1}), but expect that concentration inequalities
should be useful in this context.
\end{longlist}
\end{remarks*}
%


\subsection{Notation and scaling relations}\label{ScaRe}

Motivated by remark (ii) above, we place
the $\psi$-growth process on its proper scale by re-scaling as
follows. With $\beta$ and $\gamma$ as in (\ref{twodef}), for a
\textit{fixed} $x \in A$ and any generic point $\bar{y} := (y,h_y)
\in A_+$, we put
$\bar{y}^{(\lambda)}:= \bar{y}':= (y',h'_y)$
with
%
%
\begin{equation}\label{RESCALING}
y':= y^{({\lambda})} := \lambda^{\beta} (y-x) \quad\mbox{and}
\quad h_y':=
h_y^{({\lambda})} := \lambda^{\gamma} h_y.
\end{equation}
Also, for readability, in our notation \textit{we will
not explicitly indicate the dependency of the scaling in \textup{(\ref
{RESCALING})} on
$x.$}
The versions of $\psi, \rho,{\mathcal{P}}_{{\lambda}\rho}$ and
$\xi$ under this
re-scaling are determined by the relations
%
%
\begin{eqnarray}\label{RESCALING3}
\psi^{({\lambda})}(l) &:=& {\lambda}^{\gamma} \psi({\lambda
}^{-\beta} l),
\\
%
%
\label{rdef}
\rho^{({\lambda})} (y',h_y') &:=& {\lambda}^{\delta\gamma} \rho(y,h_y),
\\
\label{RESCALING4}
{\mathcal{P}}^{({\lambda})}_{{\lambda}\rho} &:=& {\mathcal
{P}}^{({\lambda})}_{{\lambda}\rho}[x] :=
\{ (y',h_y')\dvtx (y,h_y) \in{\mathcal{P}}_{{\lambda}\rho} \}
\end{eqnarray}
and
%
%
\begin{equation}\label{RESCALING5}
\xi^{({\lambda})} ((y',h'_y), \{ (y_i',h'_{y_i}) \}_{i \geq
1} )
:= \xi( (y,h_y),\{(y_i,h_{y_i})\}_{i \geq1} )
\end{equation}
and likewise for $\hat{\xi}.$ Since $dy' = {\lambda}^{{\beta
}(d-1)}\,dy$
and $dh'_y = {\lambda}^{\gamma}\,dh_y$, it follows that
\[
\rho^{({\lambda})}(y',h_y') \,dy'\, dh_y' = {\lambda}\rho(y,h_y) \,
dy \,dh_y.
\]
Note also that
%
%
\begin{equation}\label{PPequiv}
\mathcal{P}_{{\lambda}\rho
}^{({\lambda})} \stackrel{\mathcal{D}}{=}
{\mathcal{P}}_{\rho^{({\lambda})}}.
\end{equation}

Moreover, by (\ref{RESCALING}) and (\ref{rdef}), $
\rho^{({\lambda})}(y',h_y') (h_y')^{-\delta}= {\lambda}^{\delta
\gamma}
\rho(y,h_y)({\lambda}^{\gamma}h_y)^{-\delta}$, where $y = {\lambda
}^{-\beta}
y' + x$. Under the above re-scaling for each fixed $x \in A$ and
for each $(y', h'_y)$, we have the crucial limit
%
%
\begin{equation}\label{SCALINGH}
\lim_{{\lambda}\to\infty} \rho^{({\lambda})}(y',h_y')
(h_y')^{-\delta}= \lim_{{\lambda}\to\infty} \rho(y,h_y)
(h_y)^{-\delta}= \rho_0(x)
\end{equation}
and by ($\Psi$2) and (\ref{RESCALING3}), for all $ l \in\mathbb{R}_+$,
%
%
\begin{equation}\label{SCALINGPSI}
\ \lim_{{\lambda}\to\infty} \psi^{({\lambda})}(l) = l^{\alpha}.
\end{equation}
It is also worth noting that $\xi^{({\lambda})}$ could alternatively
be defined
by following the original definition of $\xi$ with $\psi$ replaced there
by $\psi^{({\lambda})}$; the same applies for $\hat{\xi}^{({\lambda})}$.
Observe that in fact it states approximate self-similarity of $\psi$-growth
processes under the re-scaling given by (\ref{RESCALING}) and (\ref
{RESCALING3}).
Motivated by this observation, we have already put $\psi^{(\infty
)}(l) := l^{\alpha}$ and now we define, for all $x \in A$
and for all $(y', h'_y) \in\mathbb{R}^{d-1} \times\mathbb{R}_+$,
%
%
\begin{equation}\label{SCALINGLIMIT}
\rho^{(\infty)}(y',h_y') :=
\rho_x^{(\infty)}(y',h_y') := \rho_0(x) (h_y')^{\delta}.
\end{equation}

\section{Applications}\label{s2}

We describe here applications of the main results. We limit the
discussion to the following:
\begin{longlist}
\item the number of vertices in the
convex hull of a Poisson sample, and
\item the number of maximal points in a Poisson or i.i.d.
sample,
\end{longlist}
but it should be emphasized that the techniques could potentially
be applied to a broader scope of examples. These include, for
instance, the variance asymptotics for Johnson--Mehl growth
processes \cite{Mo1} with nonlinear growth rates (see, e.g.,
Section 3.2.2 in \cite{BY2} for the description of the model
and the corresponding central limit theorem).
Also, as observed in Section 2.3 of \cite{Ba}, the case $\psi(l) =
l^2$ (paraboloids) may figure in the limit behavior of some point
processes associated with the asymptotic solutions of Burgers
equation
\[
\frac{\partial v} {\partial t} + v \frac{\partial v} {\partial x}
= \varepsilon\Delta v
\]
in the inviscous limit $\varepsilon\to0$. We will likewise not
treat this example either.

\subsection{Number of vertices in the convex hull of an i.i.d. sample}

Recall that $B_d$ denotes the unit radius ball centered at the
origin of $\mathbb{R}^d$ and let $\partial B_d$ denote its boundary. Let
$\rho\dvtx B_d \to\mathbb{R}_+$ be a continuous density on $B_d$. We shall
assume that $\rho(x) = \rho_0(x/|x|) (1-|x|)^{\delta} (1+o(1))$
for some $\delta\geq0$ and that $\rho_0 \dvtx\partial B_d \to
\mathbb{R}_+$
is continuous and bounded away from $0$. Let ${\mathcal{P}}_{{\lambda
}\rho}$ be a
Poisson point
process on $B_d$ with intensity measure ${\lambda}\rho(x)\, dx$ and let
$\operatorname{conv}({\mathcal{P}}_{{\lambda}\rho})$ be the random
polytope given by the convex
hull of ${\mathcal{P}}_{{\lambda}\rho}$. Recalling that $\mathcal{V}(
{\mathcal{P}_{{\lambda}
\rho} })$ denotes
the vertices of $\operatorname{conv}( {\mathcal{P}_{{\lambda}\rho} })$,
consider the \textit{vertex empirical
point measure}
%
%
\begin{equation}\label{vepm}
\mu_{{\lambda}{\rho}}:= \sum_{ x \in
\mathcal{V}( {\mathcal{P}_{{\lambda}\rho} }) } \delta_{x}.
\end{equation}

As will be shown in Section \ref{ApplSection}, Theorems \ref
{LLN}--\ref{CLT}
yield the following limit theory for ${\mu}_{ {\lambda}\rho}$. Let
$N(0,1)$ denote the standard normal random variable.

\begin{thm}\label{convexhullthm}
There are constants $M:=M(d,\delta)$ and $V:=V(d,\delta)$ such that
for all $f \in\mathcal{C}_b(B_d)$
%
\begin{eqnarray}\label{explim}
&&\lim_{\lambda\to\infty} \lambda^{-(d-1)/(d-1+2(1+\delta))}
\mathbb{E}[\langle f,
\mu_{\lambda\rho}\rangle]\nonumber\\[-8pt]\\[-8pt]
&&\qquad = M \int_{\partial B_d } f(s)
\rho_0^{(d-1)/(d-1+2(1+\delta))}(s) \,ds\nonumber
\end{eqnarray}
and
\begin{eqnarray}\label{varlim}
&&\lim_{\lambda\to\infty} \lambda^{-(d-1)/(d-1+2(1+\delta))}
\operatorname{Var}[\langle f,
\mu_{\lambda\rho}\rangle]\nonumber\\[-8pt]\\[-8pt]
&&\qquad = V \int_{\partial B_d} f^2(s) \rho
_0^{(d-1)/(d-1+2(1+\delta))}(s) \,ds.\nonumber
\end{eqnarray}
Moreover, the finite-dimensional distributions
$\lambda^{-(d-1)/2(d-1+2(1+\delta))}(\langle f_1,\bar\mu_{\lambda
\rho} \rangle,\break \ldots,
\langle f_k,\bar\mu_{\lambda\rho} \rangle),$
$f_i \in\mathcal{C}_b(B_d),$ of
$(\lambda^{-(d-1)/2(d-1+2(1+\delta))} \bar\mu_{\lambda\rho})$
converge as $\lambda\to\infty$ to those of
a mean zero Gaussian field with covariance kernel
\[
(f,g) \mapsto V \int_{\partial B_d} f(s) g(s) \rho
_0^{(d-1)/(d-1+2(1+\delta))}(s) \,ds,\qquad f,g \in\mathcal{C}_b(B_d).
\]
%
%
Additionally, if $\delta= 0$, then for all $f \in\mathcal{C}_b(B_d)$,
%
%
\begin{eqnarray}
\label{convrate}
&&\sup_t \biggl| P \biggl[ \frac{ \langle f,\bar\mu_{{\lambda
}\rho
}\rangle}{\sqrt{\operatorname{Var}\langle f, \bar\mu_{{\lambda
}\rho}\rangle
} } \leq t \biggr]
- P[N(0,1) \leq t] \biggr|\nonumber\\[-8pt]\\[-8pt]
&&\qquad= O \bigl({\lambda}^{-(d-1)/2(d + 1)}
(\log
{\lambda})^{3 + 2(d-1)} \bigr).\nonumber
\end{eqnarray}
\end{thm}

\begin{remarks*}
(i) Taking $f_1 \equiv1$ (and all other
$f_i \equiv0, i = 2,\ldots,k$) provides a central limit theorem
for the cardinality of $\mathcal{V}( {\mathcal{P}_{{\lambda}\rho} })$.

\smallskipamount=0pt
\begin{longlist}[(iii)]
\item[(ii)] Theorem \ref{convexhullthm} adds to the work of the
following authors: (a) Groeneboom \cite{Gr} and Cabo and
Groeneboom \cite{CG}, who prove a central limit theorem for the
cardinality of $\mathcal{V}({\mathcal{P}}_{{\lambda}\rho})$ when
$\rho
$ is uniform
and when $d = 2$, (b) Reitzner \cite{Re4} who considers the
one-dimensional central limit theorem and who establishes a rate of
convergence $O ({\lambda}^{-(d-1)/2(d + 1)} (\log{\lambda})^{2 + 2/(d+1)})$
to the normal for $\rho$ uniform (whence $\delta= 0$ in our
setting), without giving asymptotics for the limiting variance and
covariance, and (c) Vu \cite{VV}, who proves a central limit
theorem for the cardinality of $\mathcal{V}(\{X_i\}_{=1}^n)$, $X_i$
i.i.d. uniform, but who also does not consider limiting
covariances. Concerning rates, we believe that the power on the
logarithm, namely, $3 + 2(d-1)$, can be reduced to $2(d-1)$, but we
have not tried for this sharper rate.

\item[(iii)] As shown by Reitzner (Lemma 7 of \cite{Re4}), when
$\delta
=0$, the right-hand side of (\ref{varlim}) is strictly positive and
finite whenever $f$ is not identically zero.
\end{longlist}
\end{remarks*}

\subsection{Number of maximal points in an i.i.d. sample}\label{NMAX}

For all $\bar{w}:=(w,h_w)$, we define the downward cone
%
%
\begin{equation}\label{downarrow}
K^{\downarrow}[\bar{w}]:=\{(z,h_z) \in\mathbb{R}^{d-1} \times
\mathbb{R}_+\dvtx h_z
\leq h_w - \psi(|z - w|)\}.
\end{equation}
Consider $\psi(l) := l^{\alpha}$,
$\alpha\in(0,1],$ in Section \ref{Terminology}
so that $K[\mathbf{0}]
:= \{ (w_1,\ldots,w_d)\dvtx\break w_d \geq
(w_1^2+\cdots+w_{d-1}^2)^{\alpha/2} \}.$
Given a locally finite set $\mathcal{X}\subset\mathbb{R}^d$, a point
$\bar{w}
\in\mathcal{X}$
is called \textit{K-maximal} iff $\bar{w}$ does not belong to any $u
\oplus K[\mathbf{0}]$ for $u \in\mathcal{X}.$ When $\alpha\in
(0,1]$ we have the
equivalence $\bar y \in K[\bar x]$ iff $K[\bar y] \subseteq
K[\bar x]$ and $\bar x \in K^{\downarrow}[\bar y]$ iff
$K^{\downarrow}[\bar x] \subseteq K^{\downarrow}[\bar y]$. It thus
follows that for such $\psi$ the present notion of maximality is
just a rephrasing of the maximality notion as discussed in Section
\ref{GenRes}.
Indeed, we
see that $\bar{w}$ is $K$-maximal or $\psi$-extremal in $\mathcal
{X}$ iff
$\bar{w} \oplus K^{\downarrow}[\mathbf{0}]$ contains no other points in
$\mathcal{X}$. This is not the case for $\alpha> 1$, where the equivalence
$\bar y \in K[\bar x]$ iff $K[\bar y] \subseteq K[\bar x]$ does
not hold.



Recalling that $\mathcal{M}_K({\mathcal{P}_{{\lambda}\rho} })$
denotes the collection of $K$-maximal points in ${\mathcal
{P}_{{\lambda
}\rho}
}$, and with
$\rho$ and $A$ as in Section 1.1,
consider the induced \textit{maximal
point measure}
\[
\mu_{{\lambda}\rho}:= \sum_{x \in\mathcal{M_K}( {\mathcal
{P}_{{\lambda
}\rho} }) }
\delta_{x}.
\]

Recalling the definitions of $I(f)$ and $J(f)$ at (\ref{onept})
and (\ref{twoptt}), respectively, we have the following:

%
\begin{thm}\label{maximalthm}
With $\tau$ as given by \textup{(\ref{TAUU})} and $\alpha\in
(0,1]$, for all
$f \in\mathcal{C}_b(A_+)$,
%
%
\begin{equation}
\label{explimmax} \lim_{{\lambda}\to\infty} {\lambda}^{-\tau}
\mathbb{E}[\langle f,
\mu_{{\lambda}\rho}\rangle] = I(f)
\end{equation}
and
%
%
\begin{equation}\label{varlimmax}
\lim_{{\lambda}\to\infty}
{\lambda}^{-\tau} \operatorname{Var}
[\langle f,
\bar\mu_{{\lambda}\rho}\rangle]
= I(f^2) + J(f^2).
\end{equation}
Moreover, the finite-dimensional distributions
$(\langle f_1, {\lambda}^{-\tau/2} \bar\mu_{{\lambda}\rho}
\rangle, \ldots,\break
\langle f_k, {\lambda}^{-\tau/2} \bar\mu_{{\lambda}\rho}
\rangle),$
$f_1,\ldots,f_k \in\mathcal{C}_b(A_+),$ of
$ {\lambda}^{-\tau/2} \bar\mu_{{\lambda}\rho}$ converge as
${\lambda}\to\infty$
to those of
a mean zero Gaussian field with covariance kernel
\[
(f,g) \mapsto I(fg) + J(fg), \qquad f, g \in\mathcal{C}_b(A_+).
\]
Additionally, if $\delta= 0$, then for all $f \in\mathcal{C}_b(A_+)$,
%
%
\begin{eqnarray}\label{maxptrate}
&&\sup_t \biggl| P \biggl[ { \langle
f,\bar\mu_{{\lambda}\rho}\rangle\over\sqrt{\operatorname
{Var}\langle f,
\bar\mu_{{\lambda}\rho}\rangle} } \leq t \biggr] - P[N(0,1) \leq t]
\biggr| \nonumber\\[-8pt]\\[-8pt]
&&\qquad= O \bigl({\lambda}^{-(d-1)/2d} (\log{\lambda})^{3 +
2(d-1)} \bigr).\nonumber
\end{eqnarray}
\end{thm}

Theorem \ref{maximalthm} admits de-Poissonization as follows. Let
$X_1,\ldots,X_n$ be i.i.d. chosen in $A_+$ according to
the density $\rho,$ assumed to be integrable to $1,$ and
consider the associated maximal point measure
\[
\nu_{n}^{\xi} := \sum_{x \in\mathcal{M_K}(\{ X_i \}_{i=1}^n)}
\delta_x.
\]
%

We have then the following equivalent of Theorem \ref{maximalthm}
for binomial samples.
\begin{thm}\label{maximalthmBin}
With $\tau$ as given by \textup{(\ref{TAUU})} and $\alpha\in(0,1],$
for all $f \in\mathcal{C}_b(A_+)$,
%
%
\begin{equation}
\label{explimmaxBin}
\lim_{n \to\infty} n^{-\tau} \mathbb{E}[\langle f,
\nu_{n}^{\xi} \rangle] = I(f)
\end{equation}
and
%
%
\begin{equation}\label{varlimmaxBin}
\lim_{n \to\infty} n^{-\tau}\operatorname{Var}
[\langle f,\bar\nu_{n}^{\xi} \rangle]
= I(f^2) + J(f^2).
\end{equation}
Moreover, the finite-dimensional distributions
$(\langle f_1, n^{-\tau/2} \bar\nu_{n}^{\xi} \rangle, \ldots,
\langle f_k, n^{-\tau/2} \bar\nu_{n}^{\xi} \rangle),$
$f_1,\ldots,f_k \in\mathcal{C}_b(A_+),$ of
$ n^{-\tau/2} \bar\nu_{n}^{\xi}$ converge as $n \to\infty$ to
those of
a mean zero Gaussian field with covariance kernel
\[
(f,g) \mapsto I(fg) + J(fg), \qquad f, g \in\mathcal{C}_b(A_+).
\]
\end{thm}

\begin{remark*}
Theorems \ref{maximalthm} and
\ref{maximalthmBin} extend and generalize the work of (a) Barbour
and Xia \cite{BX}, who establish central limit theorems for the
case of homogeneous spatial temporal input, with $K$ the positive
octant in $\mathbb{R}^d$, and who consider neither convergence of
finite-dimensional distributions nor convergence of variances,
(b)~Baryshnikov and Yukich \cite{BY2}, who establish convergence of
finite-dimensional distributions but who restrict to homogeneous
temporal input (${\delta}= 0$) as well as to the case $\psi(l) = l$
(i.e., $\alpha= 1$), and (c) Baryshnikov \cite{Ba}, who also
restricts to homogeneous temporal input and does not consider
convergence of finite-dimensional distributions.
\end{remark*}


\section{Proof of main results}\label{s3}

In this section we prove Theorems \ref{LLN}--\ref{CLT}. An essential
component of the proofs involves
introducing a notion of \textit{ localization}, which quantifies the
decoupling
property of the considered functional $\xi$ over distant regions.
It is straightforward to check that the proofs hold for
$\psi$-growth without overlap when $\psi$ is linear.

\subsection{Stabilization for $\Psi$-growth functionals}\label{PSISTA}

With $B_{d-1}(y,r)$ standing 
as usual for the $(d -1)$-dimensional ball
centered at $y \in\mathbb{R}^{d-1}$ with radius $r \in(0, \infty)$, we
denote by $C_{d-1}(y,r)$ the cylinder $B_{d-1}(y,r) \times\mathbb{R}_+$.
Recalling $\bar{y}:= (y,h_y)$, consider for all $r > 0$ the finite
range version of $\xi(\bar{y},\mathcal{X})$, namely,
\[
\xi_{[r]}(\bar{y},\mathcal{X}):= \xi_{C_{d-1}(y,r)} (\bar
{y},\mathcal{X}),
\]
that is, $\xi_{[r]}(\bar{y},\mathcal{X})$ depends only on the local
behavior of $\mathcal{X}$ with
\textit{ spatial} coordinates restricted to the $r$-neighborhood of $y.$
For a
point process $\mathcal{P}$ (usually chosen to be Poisson in the sequel)
in $\mathbb{R}^{d-1}
\times\mathbb{R}_+,$ the \textit{localization radius} of $\xi$ at
$\bar{y}
\in\mathbb{R}^{d-1}
\times\mathbb{R}_+$ is defined by
%
%
\begin{equation}\label{RoS}
R^{\xi} := R^{\xi}[\bar{y};\mathcal{P}] := \inf\bigl\{ r
\in\mathbb{R}_+\dvtx\forall
s \geq r\ \xi(\bar{y},\mathcal{P}) =
\xi_{[s]}(\bar{y},\mathcal{P}) \bigr\}.
\end{equation}
%

In full analogy with $\xi^{({\lambda})}$ given by (\ref
{RESCALING5}), we
define for all ${\lambda}> 0$ the localization radius
$R^{\xi^{({\lambda})}}[ \cdot; \cdot]$ by
\[
R^{\xi^{({\lambda})}} := R^{\xi^{({\lambda})}} [\bar{y};\mathcal{P'}]
:= \inf\bigl\{ r \in\mathbb{R}_+\dvtx\forall s \geq r\  \xi
^{({\lambda})}(\bar
{y}',\mathcal{P}') = \xi^{({\lambda})}_{[s]}(\bar{y}',\mathcal
{P}') \bigr\}.
\]

Observe that the localization radius considered here formally
differs from the stabilization radii
considered in \cite{BY2}, \cite{Pe1,PY2,PY4,PY5}, essentially
defined for all $\bar{y}:=(y,h)$
to be the smallest positive real $r$ such that
$\xi(\bar{y}, (\mathcal{P}
\cap C_{d-1}(y,r)) \cup{\mathcal{A}}) = \xi(\bar{y}, (\mathcal{P}
\cap C_{d-1}(y,r)) $
for all finite ${\mathcal{A}}\subset C^c_{d-1}(y,s)$.
However, the $\psi$-extremal functional is in general extremely
sensitive to the choice of the ``outside'' configuration ${\mathcal
{A}}\subset
C^c_{d-1}(y,s)$, 
rendering the
existence and use of standard stabilization radii a bit difficult.
The benefit of the localization radius is that it considers only
the outside configurations involving points from $\mathcal{P}$.
However, since the localization radius shares many of the same
properties as the stabilization radii in \cite{BY2}, \cite{Pe1,PY2,PY4,PY5},
\textit{we will abuse terminology and henceforth refer to
the localization radius $R^{\xi}$ as a stabilization radius}.

The following lemma shows that $\xi^{({\lambda})}$ given by (\ref
{RESCALING5}) has a stabilization radius
whose tail decays exponentially uniformly in large enough ${\lambda}$ when
$\mathcal{P}$ is ${\mathcal{P}}_{{\lambda}
\rho}^{({\lambda})}$ given by (\ref{RESCALING4}) or when $\mathcal
{P}$ is
given by ${\mathcal{P}}_{{\lambda}\rho}^{({\lambda})}
\cup\{ \bar z'_1,\ldots,\bar z'_k\}$, $k \geq1,$ where $\bar z'_i,\;
i=1,\ldots,k$,
are certain deterministic points (fixed atoms). This result will
prove useful later in showing exponential decay of correlation
functions for $\psi$-growth processes.

\begin{lemm}\label{StabLemma}
\textup{(i)} For $A$ compact and convex,
there exists a constant $C$ such that, uniformly in $x$
and ${\lambda}$
large enough, for all $\bar{y}' \in{\lambda}^\beta A \times\mathbb{R}_+$
and for all collections $\{ \bar z'_1,\ldots, \bar z'_k \} \subseteq
{\lambda}^{\beta} A \times\mathbb{R}_+$ of deterministic points, $k
\geq0,$
we have for all $L > 0$
%
%
\begin{equation}\label{SRDecay}
P \bigl[ R^{\xi^{({\lambda})}} \bigl[\bar{y}';{\mathcal
{P}^*}_{{\lambda
}\rho}^{({\lambda})}
\bigr] > L \bigr]
\leq C\exp\biggl(- \frac{L^{\alpha+ d -1}}{C} \biggr),
\end{equation}
where ${\mathcal{P}^*}_{{\lambda}\rho}^{({\lambda})}:= \mathcal
{P}_{{\lambda}\rho}^{({\lambda}
)} \cup
\{ \bar z_1,\ldots,z_k\}$, so that, in particular,
${\mathcal{P}^*}_{{\lambda}\rho}^{({\lambda})} = \mathcal
{P}_{{\lambda
}\rho}^{({\lambda})}$ for
$k=0.$

\textup{(ii)} An identical bound holds if instead $A := \mathbb
{R}^{d-1}$ and
$\mathcal{P}_{{\lambda}\rho
}^{({\lambda})}$ is replaced by a homogeneous Poisson point process on
$\mathbb{R}^{d-1}$.
\end{lemm}

\begin{remark*}
In place of (\ref{SRDecay}) we have uniformly in
$x$ and ${\lambda}$ large enough, for all $\bar{y}' \in{\lambda
}^\beta A
\times\mathbb{R}_+$ and for all $L > 0$, the simpler bound
%
%
\begin{equation}
\label{SRDecay1}
P \bigl[ R^{\xi^{({\lambda})}} \bigl[\bar{y}';{\mathcal
{P}^*}_{{\lambda
}\rho}^{({\lambda})}
\bigr] > L \bigr]
\leq C\exp\biggl(- \frac{L}{C} \biggr).
\end{equation}
\end{remark*}

\begin{pf*}{Proof of Lemma \ref{StabLemma}}
We will only prove Lemma
\ref{StabLemma}(i) as identical arguments handle Lemma
\ref{StabLemma}(ii). Also, since the proof relies on probability
bounds for certain regions being devoid of points of the underlying
point process ${{\mathcal{P}^*}}_{{\lambda}\rho}^{({\lambda})},$ as
easily noted below,
we can assume without loss of generality that $k=0$ so that
${\mathcal{P}^*}_{{\lambda}\rho}^{({\lambda})} = \mathcal
{P}_{{\lambda
}\rho}^{({\lambda})}.$
Moreover, to simplify the argument below, we ignore the boundary
effects arising when
$\bar{y}'$ is close to $\partial({\lambda}^{\beta}A \times\mathbb
{R}_+),$ noting
that the absence of points
of $\mathcal{P}_{{\lambda}\rho}^{({\lambda})}$ in the vicinity of
$\bar
{y}'$ can
only decrease
$ R^{\xi^{({\lambda})}} [\bar{y}';\mathcal{P}_{{\lambda}\rho
}^{({\lambda})} ].$
This allows us to avoid obvious but technical separate considerations
for $\bar{y}'$ close to $\partial({\lambda}^{\beta}A \times\mathbb
{R}_+)$. Also,
we consider $x$ fixed but arbitrary, keeping in mind that the required
uniformity
in $x$ follows by the boundedness of $\rho,$ both from above and away
from $0.$

Define for fixed
$\bar{y}':= (y',h_y')$ and all ${\lambda}\in[0,\infty]$ the \textit{
scaled} upward cone
%
%
\begin{equation}\label{upconeScaled}
K^{({\lambda})}[\bar{y}']:= \bigl\{
(v',h_v') \in\mathbb{R}
^{d-1} \times\mathbb{R}_+\dvtx
h'_v \geq h'_y + \psi^{({\lambda})}(|v' - y'|) \bigr\}
\end{equation}
and the \textit{scaled} downward cone
%
%
\begin{equation}\label{downconeScaled}
K^{\downarrow}_{({\lambda})}[\bar{y}'] := \bigl\{ (v',h'_v) \in
\mathbb{R}^{d-1}
\times\mathbb{R}_+ \dvtx
h'_v \leq h'_y - \psi^{({\lambda})}(|v'-y'|) \bigr\}.
\end{equation}
Note that $\bar{u}' \in K^{({\lambda})} [\bar{z}']$ iff $h_u' \geq h_z'
+ \psi(|u' - z'|)$, which is equivalent to $h_z' \leq h_u' -
\psi(|z' - u'|)$, and thus, the \textit{duality} $\bar{u}' \in
K^{({\lambda})}[\bar{z}']$ iff $\bar{z}' \in
K^{\downarrow}_{({\lambda})}[\bar{u}']$.

To proceed, note that the event
$\{ R^{\xi^{({\lambda})}} [\bar{y}';\mathcal{P}_{{\lambda}\rho
}^{({\lambda})} ] > L \}
$ is equivalent
to the event
\[
E:= \bigl\{ \exists{r > L}\dvtx\ \xi^{({\lambda})} \bigl(\bar
{y}',\mathcal{P}_{{\lambda}
\rho}^{({\lambda})}\bigr)
\neq\xi_{[r]}^{({\lambda})}\bigl( \bar{y}',\mathcal{P}_{{\lambda
}\rho
}^{({\lambda})} \bigr) \bigr\},
\]
and moreover, $E \subset E_1 \cup E_2$, where $E_1$ and $E_2$ are
defined below.
Roughly speaking, the event $E_1$ ensures that $\bar{y}'$ is extremal with
respect to ${\mathcal{P}}_{{\lambda}\rho}^{({\lambda})} \cap
C_{d-1}(y',r)$ for some $r >
L$ but not
necessarily with
respect to ${\mathcal{P}}_{{\lambda}\rho}^{({\lambda})}$, whereas
$E_2$ is just the opposite.
\begin{description}
\item[Event $E_1$\textup{:}] For some $r > L$, there exists a
boundary point
$\bar{u}' \in\partial( K^{({\lambda})} [\bar{y}']) \cap C_{d-1}(y',r)$,
and such that $\bar{u}'\notin\bigcup_{\bar{z}'
\in[\mathcal{P}_{{\lambda}\rho}^{({\lambda})} \setminus\{ \bar
{y}'\}
] \cap C_{d-1}(y',r)}
K^{({\lambda})}[\bar{z}']$ but $\bar{u}' \in
\bigcup_{\bar{z}'
\in\mathcal{P}_{{\lambda}\rho}^{({\lambda})} \cap C_{d-1}(y',r)}
K^{({\lambda})}[\bar{z}']$, that is, $\xi_{[r]}^{({\lambda})}( \bar
{y}',\mathcal{P}_{{\lambda}\rho}^{({\lambda})}
)= 1$, but possibly\break $ \xi^{({\lambda})} (\bar{y}',\mathcal
{P}_{{\lambda
}\rho}^{({\lambda})})= 0$.
\item[Event $E_2$\textup{:}] For some $r > L$, there exists a
boundary point
$\bar{u}' \in\partial( K^{({\lambda})} [\bar{y}']) \cap C^c_{d-1}(y',r)$
such that $\bar{u}' \notin
\bigcup_{\bar{z}'
\in\mathcal{P}_{{\lambda}\rho}^{({\lambda})} \setminus\{ \bar
{y}'\}
} K^{({\lambda})}
[\bar{z}'],$
but $K^{({\lambda})}[ \bar{y}'] \cap C_{d-1} (\bar{y}',r) \subset
\bigcup_{\bar{z}'
\in[\mathcal{P}_{{\lambda}\rho}^{({\lambda})} \setminus\{ \bar
{y}'\}
] \cap C_{d-1}(y',r)}
K^{({\lambda})}[\bar{z}']$, that is, $\xi^{({\lambda})} (\bar
{y}',\mathcal{P}_{{\lambda}
\rho}^{({\lambda})})
= 1$ but\break $\xi_{[r]}^{({\lambda})}( \bar{y}',\mathcal
{P}_{{\lambda}\rho
}^{({\lambda})}
)= 0$.
\end{description}

On event $E_1$ writing $\bar{u}':= (u',h'_u),$ we easily check that
%
%
\begin{equation}\label{HU}
h'_u \geq\psi^{({\lambda})} \biggl(\frac{L}{2} \biggr).
\end{equation}
Indeed, we have:
\begin{itemize}
\item either $|u'-y'| \geq r \slash2$ or
\item$d(u',\partial B_{d-1}(y',r)) \geq r \slash2$ and, hence,
$d(u',z') \geq r \slash2$
for all
$\bar{z}' \in\mathcal{P}_{{\lambda}\rho}^{({\lambda})} \cap
C^c_{d-1}(y',r)$.
\end{itemize}
In both cases, on $E_1$, $(u',h'_u)$ falls into $K^{({\lambda})}[\bar{v}']$
for some
$\bar{v}'$ such that $|v'-u'| \geq r \slash2,$ either with $\bar{v}'
= \bar{y}'$
or $\bar{v}' \in\mathcal{P}_{{\lambda}\rho}^{({\lambda})} \cap
C^c_{d-1}(y',r)$. Consequently, recalling that $r > L$ and using the definition
of $K^{({\lambda})}[\cdot]$, we obtain (\ref{HU}) as required.

On $E_1$ we have $\bar{u}'\notin\bigcup_{\bar{z}'
\in[\mathcal{P}_{{\lambda}\rho}^{({\lambda})} \setminus\{ \bar
{y}' \}]
\cap C_{d-1}(y',r)} K^{({\lambda})}[\bar{z}']$, implying that the
downward cone $ K^{\downarrow}_{({\lambda})}[\bar{u}']$ is devoid of
points of $\mathcal{P}_{{\lambda}\rho}^{({\lambda})} \ \cap\
C_{d-1}(y',r).$ By the assumed properties of $\psi$ and $\rho$,
the 
integral of $\rho^{({\lambda})}$ over $K^{\downarrow}_{({\lambda
})}[\bar{u}']$ is
$\Omega(\operatorname{Vol}(K^{\downarrow}_{({\lambda})}[\bar
{u}']))$, which is
\begin{eqnarray}\label{volinvert}
\Omega\biggl(\int_0^{h'_u} \bigl(\bigl[\psi^{({\lambda
})}\bigr]^{-1}(h'_u-h')\bigr)^{d-1}\, dh'
\biggr)
&=&
\Omega\biggl(\int_0^{h'_u} (h_u'-h')^{(d-1)/\alpha} \, dh' \biggr)
\nonumber\\[-8pt]\\[-8pt]
&=&
\Omega\bigl( (h'_u)^{(\alpha+ d - 1)/ \alpha} \bigr),\nonumber
\end{eqnarray}
with the second equality
following by the definition of $[\psi^{({\lambda})}]^{-1}$,
and where
we use $f({\lambda}) = \Omega(g({\lambda}))$ to signify that
$f({\lambda})/g({\lambda})$ is asymptotically bounded away from zero.
Clearly,
the integral of $\rho^{({\lambda})}$ over $K^{\downarrow}_{({\lambda
})}[\bar
{u}'] \cap C_{d-1}(y',r)$
for $\bar{u}' \in C_{d-1}(y',r)$ is of
the same order.

Recalling from (\ref{PPequiv})
that the intensity measure of the Poisson process $\mathcal
{P}_{{\lambda}
\rho}^{({\lambda})}$ has its density given by $\rho^{({\lambda})}$,
we thus conclude for fixed $\bar{u}'$ that the probability of
the considered event $\Xi[\bar{u}'] :=
\{ K^{\downarrow}_{({\lambda})}[\bar{u}'] \cap[\mathcal
{P}_{{\lambda
}\rho}^{({\lambda})}
\setminus\{ \bar{y}' \}] \cap C_{d-1}(y',r)
= \varnothing\}$ satisfies
%
%
\begin{equation}\label{AddedEqn}
P[\Xi[\bar{u}']] \leq\exp\bigl( -\Omega\bigl( (h'_u)^{(\alpha+
d - 1)/
\alpha}\bigr)\bigr).
\end{equation}
To proceed, we recall that $r > L$ and we partition $\mathbb{R}^{d-1}
\times\mathbb{R}_+$ into unit volume cubes and we let $q_1,q_2,\ldots$
be an enumeration of those cubes having nonempty intersection with
$\partial(K^{({\lambda})} [\bar{y}'])$. Let
\[
p_i := P \bigl[ \exists{\bar{u}' \in q_i}\dvtx K^{\downarrow
}_{({\lambda})}[\bar{u}']
\cap\bigl[\mathcal{P}_{{\lambda}\rho}^{({\lambda})} \setminus\{
\bar{y}'
\}\bigr] \cap
C_{d-1}(y',L) = \varnothing\bigr]
\]
for all $i=1,2,\ldots$ and
note that, by (\ref{AddedEqn}), we have
\[
p_i \leq
\exp\bigl(-\Omega\bigl((h'_q)^{(\alpha+ d - 1)/\alpha}\bigr
)\bigr),
\]
where $h'_{q_i}$
is the last coordinate of the center of the cube $q_i.$

We now have
%
\[
P[E_1] \leq\sum_{i=1}^{\infty} p_i \leq C \int_{\psi^{({\lambda}
)}(L/2)}^{\infty} L^{d-2}
\exp\biggl( -\frac{1}{C} (h'_u)^{(\alpha+ d - 1)/\alpha} \biggr
)\, dh'_u
\]
for some $0 < C <\infty$ in view of the discussion above. Here
$CL^{d-2}$ bounds the number of cubes in the set $q_1,q_2,\ldots$ of
any fixed height $h'_u \geq
\psi^{({\lambda})}(L/2)$.

Recalling that $\psi^{({\lambda})}(L/2) = (1+o(1))(L/2)^{\alpha}$,
it follows (using a different choice of $C$ if
necessary) that
\[
P[E_1] \leq C \exp\biggl(- \frac{1}{C} L^{\alpha+ d - 1} \biggr).
\]

To estimate $P[E_2]$,
note that for $\bar{u}' := (u',h'_u) \in
\partial(K[\bar{y}'])$ lying in $C_{d-1}^c(y',r)$ we must have
%
%
\begin{equation}\label{HU2}
h'_u \geq\psi^{({\lambda})}(r).
\end{equation}
Further, since $\bar{u}' \notin\bigcup_{\bar{z}'
\in\mathcal{P}_{{\lambda}\rho}^{({\lambda})} \setminus\{ \bar
{y}' \}
} K^{({\lambda})}
[\bar{z}']$,
we have $K^{\downarrow}_{({\lambda})}[\bar{u}'] \cap[\mathcal
{P}_{{\lambda}\rho
}^{({\lambda})} \setminus\{ \bar{y}' \}]
= \varnothing.$ Denoting this event $\Xi^*[\bar{u}'] :=
\{ K^{\downarrow}_{({\lambda})}[\bar{u}'] \cap[\mathcal
{P}_{{\lambda
}\rho}^{({\lambda}
)} \setminus\{ \bar{y}' \}]
= \varnothing\}$, noting that as in (\ref{AddedEqn}) we have
%
%
\begin{equation}\label{AddedEqn2}
P[\Xi^*[\bar{u}']] \leq\exp\bigl( -\Omega\bigl( (h'_u)^{(\alpha
+ d - 1)/
\alpha}\bigr)\bigr),
\end{equation}
recalling that $r>L$ and proceeding in analogy with the case of
event $E_1$ above, with (\ref{HU}) and (\ref{AddedEqn}) there
replaced by (\ref{HU2}) and (\ref{AddedEqn2}) respectively and
with $C^c_{d-1}(y',L)$ partitioned into unit volume cubes, we
bound $P[E_2]$ by 
%
\[
P[E_2] \leq C \int_{s = L}^{\infty} s^{d-2} \int_{h_y' +
\psi^{({\lambda})}(s)}^{\infty} \exp\biggl( -\frac{1}{C}
(h'_u)^{(\alpha
+ d - 1)/ \alpha} \biggr)\,dh_u'\,ds
\]
for some $0 < C < \infty$.
It follows that $P[E_2] \leq C \exp(-{L^{\alpha+ d - 1}/ C}
). $ Since $ P [ R^{\xi^{({\lambda})}} [\bar{y}';\mathcal
{P}_{{\lambda}
\rho}^{({\lambda})} ] > L ] = P[E] \leq P[E_1] + P[E_2]$, Lemma
\ref{StabLemma} follows.
\end{pf*}

Given $\bar{y}:=(y',h'_y)$, we expect for
large temporal $h'_y$, that $\bar{y}$ is $\psi$-extremal with small
probability. Also, as previously noted in Section \ref{Terminology}, we
expect for small $\alpha$ that $\bar{y}$ is more likely to be
$\psi$-extremal. The next lemma makes these probabilities a bit
more precise and shows that the probability
of having $(y',h_y')$ extreme in ${\mathcal{P}^*}_{{\lambda}\rho
}^{({\lambda})}:=
\mathcal{P}_{{\lambda}\rho}^{({\lambda})} \cup\{ \bar z'_1,\ldots
,\bar z'_k\},
k\geq0,$
with respect to $\psi^{({\lambda})}$ decays exponentially with $h_y'$
uniformly in ${\lambda}$ for ${\lambda}$ large enough.

\begin{lemm}\label{expbds}
There exists a constant $C$ such that,
uniformly in ${\lambda}$ large
enough, for all $\bar{y}' \in{\lambda}^\beta A \times\mathbb{R}_+$
and $\{ \bar z'_1,\ldots,\bar z'_k \}$, we have
\[
P \bigl[ \xi^{({\lambda})}\bigl(\bar{y}',{\mathcal
{P}^*}_{{\lambda}\rho
}^{({\lambda})}\bigr) = 1
\bigr]
\leq C\exp\biggl(-\frac{1}{C}(h'_y)^{(\alpha+ d - 1)/ \alpha}
\biggr).
\]
\end{lemm}

\begin{pf}
Clearly, since adding extra points to $\mathcal{P}_{{\lambda}\rho
}^{({\lambda})}$ decreases
the probability of $(y',h_y')$ being extreme, we may without loss of generality
choose $k=0$ so that ${\mathcal{P}^*}_{{\lambda}\rho}^{({\lambda})} =
\mathcal{P}_{{\lambda}\rho}^{({\lambda})}.$

On the event
$E:=
\{ \xi^{({\lambda})}(\bar{y}',\mathcal{P}_{{\lambda}\rho
}^{({\lambda
})}) = 1 \}$ there
exists $\bar{u}':= (u',h'_u) \in
\partial(K^{({\lambda})}[\bar{y}'])$ such that $\bar{u}' \notin
\bigcup_{\bar{z}' \in\mathcal{P}_{{\lambda}\rho}^{({\lambda})}
\setminus\{ \bar
{y}' \}} K^{({\lambda})}[\bar{z}'],$
which is equivalent
to $K^{\downarrow}_{({\lambda})}[\bar{u}'] \cap
[\mathcal{P}_{{\lambda}\rho}^{({\lambda})} \setminus\{ \bar{y}' \}
] =
\varnothing
.$ As in the proof of Lemma 3.1,
for fixed $\bar{u}'$, the probability of the last event does not exceed
\[
\exp\biggl[-\int_{K^{\downarrow}_{({\lambda})}[\bar{u}']} \rho
^{({\lambda}
)}(v'h'_v)\,dv'\,dh'_v
\biggr] \leq C\exp\biggl(-\frac{1}{C}(h'_u)^{(\alpha+ d - 1)/
\alpha} \biggr).
\]
Recalling the relation $h'_u = h'_y + \psi^{({\lambda})}(|u'-y'|),$ putting
$|u'-y'| = s$, and
resorting again to a partition of $\mathbb{R}^{d-1} \times\mathbb
{R}_+$ into
unit volume cubes and summing up the respective probabilities as in the
proof of Lemma
\ref{StabLemma}, we obtain the required bound
\begin{eqnarray*}
P[E] &\leq& C \int_0^{\infty} s^{d-2} \int_{h_y'}^{\infty}
\exp\biggl( -\frac{1}{C} (h'_u)^{(\alpha+ d - 1)/ \alpha}
\biggr)\,dh_u'\,ds \\
&\leq& C\exp\biggl(-\frac{1}{C}(h'_y)^{(\alpha+ d -
1)/ \alpha} \biggr).
\end{eqnarray*}
\upqed
\end{pf}

\subsection{\texorpdfstring{Proof of Theorem \textup{\protect\ref
{LLN}}}{Proof of Theorem 1.1}}\label{LLNproof}

Recall the definition of ${\mathcal{P}}_{\rho_x^{(\infty)}}$ from
(\ref{SCALINGLIMIT}).
One benefit of stabilization is
that the one point correlation function\break
$\mathbb{E} [\xi^{(\infty)} ((\mathbf{0},h'),{\mathcal
{P}}_{\rho_x^{(\infty)}}
) ]$ is approximated for large $r$ by the finite range
version
\[
\mathbb{E} \bigl[\xi^{(\infty)}_{[r]} \bigl((\mathbf
{0},h'),{\mathcal{P}}_{\rho_x^{(\infty)}}
\bigr) \bigr]
\]
and, similarly, $\mathbb{E}
[\xi^{({\lambda})} ((\mathbf{0},h'), {\mathcal
{P}}_{{\lambda}\rho}^{({\lambda})}
) ]$ is approximated by its finite range version $\mathbb{E}
[\xi_{[r]}^{({\lambda})} ((\mathbf{0},h'), {\mathcal
{P}}_{{\lambda}\rho}^{({\lambda})}
) ]$. Using the large ${\lambda}$ weak convergence of
${\mathcal{P}}_{{\lambda}\rho}^{({\lambda})}$ to ${\mathcal
{P}}_{\rho_x^{(\infty)}}$, one may
approximate the first mentioned finite range version by the second
and thus show that $\mathbb{E} [\xi^{({\lambda})}
((\mathbf{0},h'), {\mathcal{P}}_{{\lambda}
\rho}^{({\lambda})} ) ]$ is asymptotically equal to
$\mathbb{E} [\xi^{(\infty)} ((\mathbf{0},h'),{\mathcal
{P}}_{\rho_x^{(\infty)}}
) ]$. This is spelled out in Lemma \ref{L1conv}
below, which captures the essence of stabilization and which lies
at the heart of the proof of Theorem \ref{LLN}. Note that when
Lemma~\ref{L1conv} is combined with Lemma \ref{expbds}, then it
shows
%
%
\begin{equation}\label{IFF}
\mathbb{E} \bigl[ \xi^{(\infty)}
((\mathbf{0},h'),{\mathcal{P}}^* )
\bigr] \leq C\exp\biggl(-\frac{1}{C}(h')^{(\alpha+ d - 1)/
\alpha} \biggr)
\end{equation}
and, therefore, $I(f) < \infty$ for $f \in\mathcal{C}_b(A_+)$. Recall
from (\ref{RESCALING5}) that $\xi^{({\lambda})}$ is the re-scaled
version of $\xi$ with dependency on $x$ fixed.

\begin{lemm}\label{L1conv}
For all $x \in A$ and $h' \in\mathbb
{R}_+$, we
have
\[
\lim_{{\lambda}\to\infty} \mathbb{E} \bigl[\xi^{({\lambda
})} \bigl((\mathbf{0},h'),
{\mathcal{P}}_{{\lambda}\rho}^{({\lambda})} \bigr) \bigr]=
\mathbb{E} \bigl[
\xi^{(\infty)} \bigl((\mathbf{0},h'),{\mathcal{P}}_{\rho
_x^{(\infty)}} \bigr)
\bigr].
\]
\end{lemm}

\begin{pf}
Fix $x \in A$. Taking into account (\ref{SCALINGH})
and (\ref{SCALINGLIMIT})
and using the results of Section 3.5 in \cite{RES} [see Proposition 3.22
or Proposition 3.19 there combined with Proposition 3.6(ii) ibidem],
we observe that as ${\lambda}\to\infty$, ${\mathcal{P}}_{{\lambda
}\rho}^{({\lambda})}$ converges
weakly to
${\mathcal{P}}_{\rho_x^{(\infty)}}$ as a point process; see ibidem.
Using Theorem 5.5 in \cite{BILL} with $h_{{\lambda}} := \xi
_{[r]}^{({\lambda}
)}((\mathbf{0},h'),\cdot)$
and $h := \xi_{[r]}^{(\infty)}((\mathbf{0},h'),\cdot)$ there, we
easily see that,
by Lemmas \ref{StabLemma} and \ref{expbds}, under the law of the limit
process ${\mathcal{P}}_{\rho_x^{(\infty)}}$,
the discontinuity event $E$ ibidem [an infinitesimal move of the point
configuration
alters the $\xi$-value for $(\mathbf{0},h)$] is contained up to an
event of
probability $0$
in the set of point configurations $\mathcal{X}$ such that either the spatial
coordinates
of two points in $\mathcal{X}$ coincide or such that there are at
least two
points $\bar{y}',\bar{y}'' \in\mathcal{X}$
such that the boundaries of the upward cones $K^{(\infty)}[\bar{y}']$
and
$K^{(\infty)}[\bar{y}'']$ [recall
(\ref{upconeScaled})] intersect
in a point lying on the boundary of the upward cone $K^{(\infty
)}[(\mathbf{0},h')],$
which clearly happens with probability $0$ under the law of
${\mathcal{P}}_{\rho_x^{(\infty)}}.$
Indeed, Lemma \ref{expbds} states that no effects coming from $h\to
\infty$
arise (no infinite range dependencies in $h$).
A similar statement in space is provided by Lemma \ref{StabLemma}.
Combining both these statements allows us to draw conclusions
from the weak convergence of point processes as we do in the above
argument; see ibidem in \cite{RES}.
Thus, Theorem 5.5 in \cite{BILL} yields
%
%
\begin{equation}\label{limit1}
\lim_{{\lambda}\to\infty} \mathbb{E} \bigl[\xi_{[r]}^{({\lambda
})} \bigl((\mathbf{0},h'), {\mathcal{P}}
_{{\lambda}\rho}^{({\lambda})}
\bigr) \bigr] = \mathbb{E} \bigl[ \xi_{[r]}^{(\infty)}
\bigl((\mathbf{0},h'),
{\mathcal{P}}_{\rho_x^{(\infty)}} \bigr) \bigr].
\end{equation}

Let $R^{\xi}:= R^{\xi^{({\lambda})}}[(\mathbf{0},h'); {\mathcal
{P}}_{{\lambda}\rho}^{({\lambda})} ].$
We have for all $r > 0$ and all ${\lambda}> 0$
\begin{eqnarray*}
&&\mathbb{E} \bigl[\xi^{({\lambda})} \bigl((\mathbf{0},h'),
{\mathcal
{P}}_{{\lambda}\rho}^{({\lambda})} \bigr)
\bigr]\\
& &\qquad= \mathbb{E} \bigl[ \xi^{({\lambda})} \bigl((\mathbf
{0},h'), {\mathcal{P}}_{{\lambda}
\rho}^{({\lambda})} \bigr) \mathbf{1}_{R^{\xi} \leq r } \bigr] +
\mathbb{E}
\bigl[\xi^{({\lambda})} \bigl((\mathbf{0},h'), {\mathcal
{P}}_{{\lambda}\rho}^{({\lambda})} \bigr)
\mathbf{1}_{R^{\xi}
> r } \bigr]
\\
&&\qquad= \mathbb{E} \bigl[ \xi_{[r]}^{({\lambda})} \bigl(
(\mathbf{0},h'),
{\mathcal{P}}_{{\lambda}\rho}^{({\lambda})}
\bigr) {\bf1}_{R^{\xi} \leq r } \bigr] + \mathbb{E}
\bigl[\xi^{({\lambda})} \bigl((\mathbf{0},h'), {\mathcal
{P}}_{{\lambda}\rho}^{({\lambda})} \bigr)
{\bf
1}_{R^{\xi}
> r } \bigr].
\end{eqnarray*}
%
By Lemma \ref{StabLemma}(i) [recall the bound (\ref{SRDecay1})],
Cauchy--Schwarz, and the boundedness of $\xi^{({\lambda})}_{[r]}$,
uniformly in large ${\lambda}$ and all $r > 0$,
\[
\mathbb{E} \bigl[\xi_{[r]}^{({\lambda})} \bigl( (\mathbf{0},h'),
{\mathcal{P}}_{{\lambda}\rho}^{({\lambda})}
\bigr) {\bf1}_{R^{\xi} > r } \bigr] \leq C\exp\biggl(- \frac
{r}{C} \biggr)
\]
for some $C$ not depending on $x.$ Likewise, uniformly in large
${\lambda},$ we have $\mathbb{E} [\xi^{({\lambda})} (
(\mathbf{0},h'), {\mathcal{P}}_{{\lambda}
\rho}^{({\lambda})} ) {\bf1}_{R^{\xi} > r } ] \leq C
\exp(-r/C )$. It follows that, for large ${\lambda}
> 0$ and all $r
> 0$,
%
%
\begin{equation}\label{DOD0}
\bigl| \mathbb{E} \bigl[\xi^{({\lambda
})} \bigl( (\mathbf{0},h'), {\mathcal{P}}
_{{\lambda}
\rho}^{({\lambda})} \bigr) \bigr] - \mathbb{E} \bigl[\xi
_{[r]}^{({\lambda})} \bigl(
(\mathbf{0},h'), {\mathcal{P}}_{{\lambda}\rho}^{({\lambda})}
\bigr) \bigr] \bigr| \leq2C
\exp\biggl(-\frac{r}{C} \biggr).
\end{equation}
%
Similarly, Lemma \ref{StabLemma}(ii) gives, for all $r > 0$,
\[
\bigl| \mathbb{E} \bigl[\xi^{(\infty)} \bigl( (\mathbf{0},h'),
{\mathcal{P}}_{\rho
_x^{(\infty)
}} \bigr) \bigr] - \mathbb{E} \bigl[\xi_{[r]}^{(\infty)} \bigl(
(\mathbf{0},h'),
{\mathcal{P}}_{\rho_x^{(\infty)}} \bigr) \bigr] \bigr| \leq2C
\exp\biggl(-\frac{r}{C} \biggr).
\]
Write
\begin{eqnarray}\label{tri-ineq}
&& \bigl| \mathbb{E} \bigl[\xi^{({\lambda})} \bigl((\mathbf{0},h'),
{\mathcal{P}}_{{\lambda}\rho}^{({\lambda})}
\bigr) \bigr] - \mathbb{E} \bigl[ \xi^{(\infty)} \bigl((\mathbf{0},h'),
{\mathcal{P}}_{\rho_x^{(\infty)}} \bigr) \bigr] \bigr|\nonumber
\\
&&\qquad\leq\bigl| \mathbb{E}
\bigl[\xi^{({\lambda})} \bigl((\mathbf{0},h'), {\mathcal
{P}}_{{\lambda}\rho}^{({\lambda})}
\bigr) \bigr] - \mathbb{E} \bigl[\xi_{[r]}^{({\lambda})}
\bigl((\mathbf{0},h'), {\mathcal{P}}
_{{\lambda}
\rho}^{({\lambda})} \bigr) \bigr] \bigr|\nonumber
\\[-8pt]\\[-8pt]
&&\quad\qquad{}+ \bigl| \mathbb{E} \bigl[\xi_{[r]}^{({\lambda})}
\bigl((\mathbf
{0},h'), {\mathcal{P}}_{{\lambda}
\rho}^{({\lambda})} \bigr) \bigr] - \mathbb{E} \bigl[ \xi
_{[r]}^{(\infty)}
\bigl((\mathbf{0},h'), {\mathcal{P}}_{\rho_x^{(\infty)}} \bigr)
\bigr] \bigr|\nonumber
\\
&&\quad\qquad{}+ \bigl| \mathbb{E} \bigl[ \xi_{[r]}^{(\infty)}
\bigl((\mathbf{0},h'),
{\mathcal{P}}_{\rho_x^{(\infty)}} \bigr) \bigr] - \mathbb
{E} \bigl[ \xi^{(\infty)}
\bigl((\mathbf{0},h'), {\mathcal{P}}_{\rho_x^{(\infty)}} \bigr)
\bigr] \bigr|.\nonumber
\end{eqnarray}

For fixed $r$, the second term on the right-hand side of
(\ref{tri-ineq}) goes to zero as ${\lambda}\to\infty$ by
(\ref{limit1}). The first and third terms are bounded above by $
2C \exp(-r/C)$. Letting $r \to\infty$ completes the proof of
Lemma \ref{L1conv}.
\end{pf}

Given Lemmas \ref{expbds} and \ref{L1conv}, we now prove Theorem
\ref{LLN} as follows. We have
\[
\mathbb{E}[ \langle f, \mu_{{\lambda}{\rho}}^{\xi} \rangle] =
\int_A
\int_0^{\infty}
f(x,h_x) \mathbb{E}[ \xi((x,h_x), {\mathcal{P}}_{{\lambda}{\rho}})
] {\lambda}\rho(x,h_x)\, dh_x \,dx.
\]

By (\ref{RESCALING5}), we have $\xi((x,h_x), {\mathcal{P}}_{{\lambda
}{\rho}}) = \xi^{({\lambda})} ((\mathbf{0},h'_x),
{\mathcal{P}}_{{\lambda}{\rho}}^{({\lambda})} )$ and by (\ref
{rdef}), we have $\rho(x,h_x)
= {\lambda}^{-\gamma\delta} \rho^{({\lambda})}(\mathbf{0},h'_x)$.
Thus, putting
$h'_x:= {\lambda}^{\gamma}h_x$ and recalling $1 - \gamma(\delta+ 1) =
\tau$ [see (\ref{TAUU}) and (\ref{twodef})], we obtain
\[
\mathbb{E}[ \langle f, \mu_{{\lambda}{\rho}}^{\xi} \rangle] =
\int_A
\int_0^{\infty} f(x,h'_x \lambda^{-\gamma})
\mathbb{E}\bigl[ \xi^{({\lambda})} \bigl((\mathbf{0},h'_x),
{\mathcal
{P}}_{{\lambda}{\rho}}^{({\lambda})} \bigr) \bigr]
{\lambda}^{\tau} \rho^{({\lambda})}(\mathbf{0},h'_x) \,dh'_x\, dx
\]
or, simply,
\[
{\lambda}^{-\tau}\mathbb{E}[ \langle f, \mu_{{\lambda}{\rho
}}^{\xi} \rangle] = \int_A
\int_0^{\infty} f(x,h'_x \lambda^{-\gamma}) \mathbb{E}\bigl[ \xi
^{({\lambda})}
\bigl((\mathbf{0},h'_x),
{\mathcal{P}}_{{\lambda}{\rho}}^{({\lambda})} \bigr) \bigr] \rho
^{({\lambda
})}(\mathbf{0},h'_x) \,dh'_x \,dx.
\]
We put
\[
g_{{\lambda}}(x,h'_x):= \mathbb{E}\bigl[ \xi^{({\lambda})} \bigl
((\mathbf
{0},h'_x), {\mathcal{P}}_{{\lambda}{\rho
}}^{({\lambda})}
\bigr) \bigr]
\rho^{({\lambda})}(\mathbf{0},h'_x).
\]
For all $x \in A$ and $h'_x \in\mathbb{R}_+$, we have by Lemma
\ref{L1conv} and (\ref{SCALINGH})
\[
\lim_{{\lambda}\to\infty} g_{{\lambda}}(x,h'_x) = \mathbb{E}\bigl
[ \xi
^{(\infty)}
\bigl((\mathbf{0},h'_x), {\mathcal{P}}_{ \rho_x^{(\infty) } }
\bigr) \bigr] \rho_0(x)
{h'}_x^{\delta}
\]
and moreover, by Lemma \ref{expbds}
for all $(x,h) \in A_+$, $ g_{{\lambda}}(x,h'_x)$ is bounded uniformly in
${\lambda}$ by the
function $(x,h') \mapsto C' (h')^{\delta} \exp(-h'/C)$, which is
integrable on $A_+$. Consequently, the dominated convergence
theorem yields
\begin{eqnarray}\label{DomConv}
&&\quad\qquad\lim_{{\lambda}\to\infty} {\lambda}^{-\tau} \mathbb
{E}[\langle
f, \mu^{\xi}_{{\lambda}
\rho} \rangle]\nonumber\\[-8pt]\\[-8pt]
&&\qquad\qquad= \int_{A} \int_{0}^{\infty} f(x,0)
\mathbb{E} \bigl[\xi^{(\infty)} \bigl((\mathbf{0},h'),{\mathcal
{P}}_{\rho_x^{(\infty)}} \cup
\{ (\mathbf{0},h') \} \bigr) \bigr] \rho_0(x)
(h')^{\delta} \,dh' \,dx.\nonumber
\end{eqnarray}
Using the scaling relations (\ref{RESCALING}), (\ref{RESCALING3}),
(\ref{SCALINGLIMIT1}) and (\ref{SCALINGLIMIT}), we see that
\begin{eqnarray}\label{SCALINGLIMIT2}
&&\xi^{(\infty)} \bigl((\mathbf{0},h'),{\mathcal{P}}_{\rho
_x^{(\infty)}} \cup\{ (\mathbf{0}
,h') \} \bigr)\nonumber\\[-8pt]\\[-8pt]
&&\qquad\stackrel{\mathcal{D}}{=}
\xi^{(\infty)} \bigl((\mathbf{0},[\rho_0(x)]^{\gamma}
h'),{\mathcal{P}}_{*} \cup\{
(\mathbf{0},[\rho_0(x)]^{\gamma} h') \} \bigr),\nonumber
\end{eqnarray}
with $\stackrel{\mathcal{D}}{=}$ standing for equality in law.
Theorem \ref{LLN} follows
by using (\ref{SCALINGLIMIT2}), changing variables $h'' := [\rho
_0(x)]^{\gamma}h'$
 in the integral in (\ref{DomConv}) and recalling that $\tau= 1 -
\gamma(\delta+ 1).$

\subsection{\texorpdfstring{Proof of Theorem \textup{\protect\ref
{VAR}}}{Proof of Theorem 1.2}}\label{VARproof}

Fix $x \in A$ and recall from (\ref{RESCALING4}) that ${\mathcal
{P}}_{{\lambda}{\rho
}}^{({\lambda})}:= {\mathcal{P}}_{{\lambda}
{\rho}}^{({\lambda})}[x]$. For all ${\lambda}
> 0$, $h' \in\mathbb{R}_+$, and $(y',h_y') \in{\lambda}^{\beta}A
\times\mathbb{R}_+$,
consider the pair correlation function for the re-scaled growth
process:
%
\begin{eqnarray}\label{twopt}
&&c^{({\lambda})}((\mathbf{0},h'), (y',h_y'))\nonumber\\
&&\qquad:=
c_x^{({\lambda})}((\mathbf{0},h'),(y',h_y'))\nonumber\\[-8pt]\\[-8pt]
&&\qquad:=
\mathbb{E} \bigl[
\xi^{({\lambda})} \bigl((\mathbf{0},h'), {\mathcal{P}}_{{\lambda
}{\rho}}^{({\lambda})} \cup
(y',h_y') \bigr) \xi^{({\lambda})} \bigl((y',h_y'), {\mathcal
{P}}_{{\lambda}{\rho
}}^{({\lambda})}
\cup(\mathbf{0},h') \bigr) \bigr]\nonumber\\
&&\quad\qquad{}-
\mathbb{E}\xi^{({\lambda})} \bigl((\mathbf{0},h'),
{\mathcal{P}}_{{\lambda}{\rho}}^{({\lambda})} \bigr) \mathbb
{E}\xi^{({\lambda})} \bigl((y',h_y'),
{\mathcal{P}}_{{\lambda}
{\rho}}^{({\lambda})} \bigr).\nonumber
\end{eqnarray}
%
Consider also the pair correlation function for the limit growth
process $\xi^{(\infty)}$:
\begin{eqnarray*}
&& c_x^{(\infty)}((\mathbf{0},h'),
(y',h_y'))\\
&&\qquad:=
\mathbb{E} \bigl[ \xi^{(\infty)} \bigl((\mathbf{0},h), {\mathcal
{P}}_{\rho_x^{(\infty) }}
\cup
(y',h_y') \bigr) \xi^{(\infty)} \bigl((y',h_y'),
{\mathcal{P}}_{\rho_x^{(\infty) }} \cup(\mathbf{0},h') \bigr)
\bigr]\\
&&\qquad\quad{}-
\mathbb{E}
\xi^{(\infty)} \bigl((\mathbf{0},h'), {\mathcal{P}}_{\rho
_x^{(\infty) }} \bigr) \mathbb{E}
\xi^{(\infty)} \bigl((y',h_y'), {\mathcal{P}}_{\rho_x^{(\infty)
}} \bigr).
\end{eqnarray*}

A second benefit of stabilization, as shown by the next lemma, is
that it facilitates convergence of pair correlation functions and
thus leads to variance asymptotics. The next lemma is the second-order
counterpart to Lemma \ref{L1conv}.

\begin{lemm}[(Convergence of two point correlation function)]\label{corlimit}
For all $(x,\break h_x):= (x,h) \in A_+$, and $(y',h'_y) \in{\lambda
}^{\beta}A
\times\mathbb{R}_+$, we have
\[
\lim_{{\lambda}\to\infty} c^{({\lambda})}_x((\mathbf{0},h'),
(y',h_y')) =
c_x^{(\infty)}((\mathbf{0},h'), (y',h_y')).
\]
\end{lemm}

\begin{pf}
In view of Lemma
\ref{L1conv}, it will suffice to show
\begin{eqnarray}\label{covlimit}
&&\quad\qquad\lim_{{\lambda}\to\infty} \mathbb{E} \bigl[ \xi
^{({\lambda
})} \bigl((\mathbf{0},h'), {\mathcal{P}}
_{{\lambda}
{\rho}}^{({\lambda})} \cup(y',h_y') \bigr) \xi^{({\lambda})}
\bigl((y',h_y'),
{\mathcal{P}}_{{\lambda}{\rho}}^{({\lambda})} \cup(\mathbf
{0},h') \bigr) \bigr]
\nonumber\\[-8pt]\\[-8pt]
&&\qquad\quad\qquad= \mathbb{E} \bigl[ \xi^{(\infty)} \bigl
((\mathbf{0},h'), {\mathcal
{P}}_{\rho_x^{(\infty) }}
\cup(y',h_y') \bigr) \xi^{(\infty)} \bigl((y',h_y'),
{\mathcal{P}}_{\rho_x^{(\infty) }} \cup(\mathbf{0},h') \bigr)
\bigr].\nonumber
\end{eqnarray}

Let
$R^{\xi}:= R^{\xi^{({\lambda})}} [(\mathbf{0},h'); {\mathcal
{P}}_{{\lambda}{\rho}}^{({\lambda})} \cup
(y',h'_y)]$ and let $R_{y'}^{\xi}:= R^{\xi^{({\lambda})}}[(y',h'_y);
{\mathcal{P}}_{{\lambda}{\rho}}^{({\lambda})} \cup(\mathbf
{0},h')].$ For all $r > 0$, we let $E_r:=
\{R_{y'}^{\xi} \leq r, \ R^{\xi} \leq r \}.$ We split the left-hand
side of (\ref{covlimit}) as
\begin{eqnarray*}
&&\mathbb{E} \bigl[ \xi^{({\lambda})}\bigl((\mathbf{0},h'),
{\mathcal
{P}}_{{\lambda}{\rho}}^{({\lambda})} \cup(y',h_y')\bigr)
\xi^{({\lambda})}\bigl((y',h_y'), {\mathcal{P}}_{{\lambda}{\rho
}}^{({\lambda})} \cup(\mathbf{0},h')\bigr) \mathbf{1}_{E_r} \bigr
] \\
&&\qquad{} +
\mathbb{E} \bigl[ \xi^{({\lambda})}\bigl((\mathbf{0},h'),
{\mathcal
{P}}_{{\lambda}{\rho}}^{({\lambda})} \cup(y',h_y')\bigr)
\xi^{({\lambda})}\bigl((y',h_y'), {\mathcal{P}}_{{\lambda}{\rho
}}^{({\lambda})} \cup(\mathbf{0},h')\bigr) \mathbf{1}_{E_r^c}
\bigr].
\end{eqnarray*}
The second expectation is bounded by $C \exp(-r/C)$ for some $C$
not depending on $x$ by Lemma \ref{StabLemma}(i)
and by Cauchy--Schwarz. By the definition of the stabilization
radius, the first is simply
\[
\mathbb{E} \bigl[ \xi_{[r]}^{({\lambda})}\bigl((\mathbf{0},h'),
{\mathcal{P}}_{{\lambda}{\rho}}^{({\lambda})} \cup
(y',h_y')\bigr)
\xi_{[r]}^{({\lambda})}\bigl((y',h_y'), {\mathcal{P}}_{{\lambda
}{\rho
}}^{({\lambda})} \cup(\mathbf{0}
,h')\bigr)\mathbf{1}_{E_r}
\bigr].
\]
Again, by Lemma \ref{StabLemma}(i) and by Cauchy--Schwarz, for all
$r > 0$, this is within $C \exp(-r/C)$ of
\[
\mathbb{E} \bigl[ \xi_{[r]}^{({\lambda})}\bigl((\mathbf{0},h'),
{\mathcal{P}}_{{\lambda}{\rho}}^{({\lambda})} \cup
(y',h_y')\bigr)
\xi_{[r]}^{({\lambda})}\bigl((y',h_y'), {\mathcal{P}}_{{\lambda
}{\rho
}}^{({\lambda})} \cup(\mathbf{0},h')\bigr)
\bigr],
\]
hence,
\begin{eqnarray}\label{DOD1}
&& \bigl| \mathbb{E} \bigl[ \xi^{({\lambda})}\bigl((\mathbf{0},h'),
{\mathcal{P}}_{{\lambda}{\rho}}^{({\lambda})}
\cup(y',h_y')\bigr)
\xi^{({\lambda})}\bigl((y',h_y'), {\mathcal{P}}_{{\lambda}{\rho
}}^{({\lambda})} \cup(\mathbf{0},h')\bigr)
\bigr]\nonumber\\
&&\qquad{} -
\mathbb{E} \bigl[ \xi_{[r]}^{({\lambda})}\bigl((\mathbf{0},h'),
{\mathcal{P}}_{{\lambda}{\rho}}^{({\lambda}
)}\cup(y',h_y') \bigr)
\xi_{[r]}^{({\lambda})}\bigl((y',h_y'), {\mathcal{P}}_{{\lambda
}{\rho
}}^{({\lambda})} \cup(\mathbf{0},h')\bigr)
\bigr] \bigr|\\
&&\qquad\leq
2 C \exp\biggl(\frac{-r}{C} \biggr)\nonumber
\end{eqnarray}
uniformly in $x.$
Now, in analogy with (\ref{limit1}), we have
\begin{eqnarray}\label{DOD3}
&&\qquad\lim_{{\lambda}\to\infty}
\mathbb{E} \bigl[ \xi_{[r]}^{({\lambda})}\bigl((\mathbf{0},h'),
{\mathcal{P}}_{{\lambda}{\rho}}^{({\lambda})} \cup
(y',h_y')\bigr)
\xi_{[r]}^{({\lambda})}\bigl((y',h_y'), {\mathcal{P}}_{{\lambda
}{\rho
}}^{({\lambda})} \cup(\mathbf{0},h')\bigr)
\bigr]\nonumber
\\[-8pt]\\[-8pt]
&&\qquad\qquad=
\mathbb{E} \bigl[ \xi_{[r]}^{(\infty)}\bigl((\mathbf{0},h'),
{\mathcal
{P}}_{\rho_x^{(\infty)}}
\cup
(y',h_y')\bigr)
\xi_{[r]}^{(\infty)}\bigl((y',h_y'), {\mathcal{P}}_{\rho
_x^{(\infty)}}
\cup(\mathbf{0}
,h')\bigr) \bigr].\nonumber
\end{eqnarray}
By Lemma \ref{StabLemma}(ii), we have for all $r > 0$
\begin{eqnarray}\label{DOD2}
&& \bigl| \mathbb{E} \bigl[ \xi^{(\infty)}\bigl((\mathbf{0},h'),
{\mathcal
{P}}_{\rho_x^{(\infty)}}
\cup(y',h_y')\bigr)
\xi^{(\infty)}\bigl((y',h_y'), {\mathcal{P}}_{\rho_x^{(\infty)}}
\cup
(\mathbf{0},h')\bigr)
\bigr]
\nonumber\\
&&\qquad{}-
\mathbb{E} \bigl[ \xi_{[r]}^{(\infty)}\bigl((\mathbf{0},h'),
{\mathcal
{P}}_{\rho_x^{(\infty)}}
\cup(y',h_y')\bigr)
\xi_{[r]}^{(\infty)}\bigl((y',h_y'), {\mathcal{P}}_{\rho
_x^{(\infty)}}
\cup(\mathbf{0},h')\bigr)
\bigr] \bigr|\\
&&\qquad\leq
2 C \exp\biggl(-\frac{r}{ C} \biggr)\nonumber
\end{eqnarray}
as in (\ref{DOD1}).
Again, note that $C$ does not depend on $x$
since $\rho_0(x)$ is bounded away from zero.
Combining (\ref{DOD1}), (\ref{DOD3}) and (\ref{DOD2}) yields
\begin{eqnarray*}
&&\limsup_{{\lambda}\to\infty} \bigl|\mathbb{E} \bigl[
\xi^{({\lambda})}\bigl((\mathbf{0},h'), {\mathcal{P}}_{{\lambda
}{\rho
}}^{({\lambda})} \cup(y',h_y')\bigr)
\xi^{({\lambda})}\bigl((y',h_y'), {\mathcal{P}}_{{\lambda}{\rho
}}^{({\lambda})} \cup(\mathbf{0},h')\bigr)
\bigr]
\\
&&\qquad{}-
\mathbb{E} \bigl[ \xi^{(\infty)}\bigl((\mathbf{0},h'), {\mathcal
{P}}_{\rho_x^{(\infty)}} \cup
(y',h_y') \bigr)
\xi^{(\infty)}\bigl((y',h_y'), {\mathcal{P}}_{\rho_x^{(\infty)}}
\cup
(\mathbf{0},h')\bigr)
\bigr] \bigr|
\\
&&\qquad\leq4 C \exp\biggl(-\frac{r}{C} \biggr)
\end{eqnarray*}
for all $r > 0.$
We conclude the proof of Lemma
\ref{corlimit} by letting $r\to\infty.$
\end{pf}

Lemma \ref{corlimit} is not enough to establish second-order
asymptotics. We will also need that $c_x^{({\lambda})}$ is bounded by
an integrable function on $A_+ \times{\lambda}^{\beta}A \times
\mathbb{R}_+$,
that is, we will need to establish the exponential
decay of the correlation function (\ref{twopt}). 
This is done in the following lemma, which combined with Lemma
\ref{corlimit}, shows that
%
%
\begin{equation}\label{JFF}
\bigl|c^{(\infty)}_x((\mathbf{0},h'), (y',h_y'))\bigr| \leq C \exp
\biggl(- \frac{1
}{C} \max\biggl( \frac{|y'|}{2}, h'_y, h' \biggr) \biggr)
\end{equation}
and,
therefore, $J(f) < \infty$ for all $f \in\mathcal{C}_b(A_+)$.

\begin{lemm} \label{corbds}
There exists a constant $C$ such that, for all ${\lambda}> 0$,
$(x,h_x):= (x,h) \in A_+$, and $(y',h'_y) \in{\lambda}^{\beta}A
\times\mathbb{R}_+$, we have
\[
\bigl|c^{({\lambda})}_x((\mathbf{0},h'), (y',h_y'))\bigr| \leq C
\exp\biggl(-
\frac{1}{C}
\max\biggl( \frac{|y'|}{2}, h'_y, h' \biggr)
\biggr).
\]
\end{lemm}

\begin{pf}
Let $r \leq|y'|/2$ and note that, by definition of $\xi
_{[r]}^{({\lambda}
)}$, we have
\begin{eqnarray*}
&&\mathbb{E} \bigl[ \xi_{[r]}^{({\lambda})}\bigl((\mathbf{0},h'),
{\mathcal{P}}_{{\lambda}{\rho}}^{({\lambda})} \cup
(y',h_y')\bigr)
\xi_{[r]}^{({\lambda})}\bigl((y',h_y'), {\mathcal{P}}_{{\lambda
}{\rho
}}^{({\lambda})} \cup(\mathbf{0},h')\bigr)
\bigr] \\
&&\qquad=
\mathbb{E} \bigl[ \xi_{[r]}^{({\lambda})}\bigl((\mathbf{0},h'),
{\mathcal{P}}_{{\lambda}{\rho}}^{({\lambda})}
\bigr) \bigr]
\mathbb{E} \bigl[ \xi_{[r]}^{({\lambda})}\bigl((y',h_y'),
{\mathcal
{P}}_{{\lambda}{\rho}}^{({\lambda}
)}\bigr) \bigr].
\end{eqnarray*}
Recalling (\ref{DOD0})
and (\ref{DOD1}), we see that
\[
\bigl|c^{({\lambda})}_x((\mathbf{0},h'), (y',h_y'))\bigr| \leq4C
\exp
\biggl(-\frac{r}{C} \biggr)
\]
for all $r \leq|y'|/2.$ In other words, putting $r = |y'|/2$ yields
for all $(x,h) \in
A_+$ and $(y',h'_y) \in{\lambda}^{\beta}A \times\mathbb{R}_+$
\[
\bigl|c^{({\lambda})}_x((\mathbf{0},h'), (y',h_y'))\bigr| \leq C
\exp\biggl(-
\frac{ |y'| }{
2C} \biggr).
\]
Appealing to Lemma \ref{expbds} shows
\[
\bigl|c^{({\lambda})}_x((\mathbf{0},h'), (y',h_y'))\bigr| \leq2 C
\exp
\biggl(-\frac{1}{C}
\max( h'_y, h')
\biggr).
\]
Combining the previous two displays
concludes the proof of Lemma \ref{corbds}.
\end{pf}

Given Lemmas \ref{corlimit} and \ref{corbds}, we now prove Theorem
\ref{VAR} as follows. 
By the Palm
theory for Poisson processes (see, e.g., Theorem 1.6 of \cite{Pe}),
we express $\operatorname{Var}[ \langle f, \mu_{{\lambda}\rho}^{\xi
} \rangle]$ as
%
%
\begin{eqnarray}\label{VarEqn}
&&{\lambda}\int_{A+} f^2(\bar{x}) \mathbb{E}[\xi(\bar{x},
{\mathcal
{P}}_{{\lambda}{\rho}} )]
{\rho}(\bar{x}) \,d\bar{x}\nonumber\\[-8pt]\\[-8pt]
&&\qquad{}+
{\lambda}^2 \int_{A_+} \int_{A_+} f(\bar{x})f(\bar{y})
c^{(1)}_x\bigl((\mathbf{0},h_x),(y-x,h_y)\bigr) {\rho}(\bar{x})
{\rho}(\bar
{y}) \,d\bar{x} \,d\bar{y},\nonumber
\end{eqnarray}
where $\bar{x}:= (x,h_x)$ and $\bar{y} := (y,h_y).$


Following verbatim the proof of Theorem \ref{LLN} shows that after
normalization by ${\lambda}^{\tau}$, the first integral converges as
${\lambda}\to\infty$ to
\[
\int_{A} \int_0^{\infty} f^2(x,0) \mathbb{E}\bigl[\xi^{(\infty
)}((\mathbf{0},h_x'),
{\mathcal{P}}_{\rho_x^{(\infty) }} ) \bigr] {\rho}_0(x)
(h_x')^{\delta}
\,dh_x' \,dx,
\]
which by the definition of $m^{(\infty)}$ and the scaling
relation (\ref{SCALINGLIMIT2}) equals
%
\[
\int_{A} \int_0^{\infty} f^2(x,0) m^{(\infty)}(\mathbf{0},h_x')
{\rho}_0^{\tau}(x) (h_x')^{\delta} \,dh_x' \,dx.
\]

Making again the usual
substitutions $y' = {\lambda}^{\beta}(y - x)$, $ h'_x = {\lambda
}^{\gamma}
h_x,$ and $h'_y = {\lambda}^{\gamma} h_y$ and recalling ${\rho}(x,h_x)=
{\lambda}^{-\gamma\delta} {\rho}^{({\lambda})}(\mathbf{0},h'_x)$,
the second integral in
(\ref{VarEqn}) becomes
\begin{eqnarray*}
&&{\lambda}^{2 - 2\gamma- 2\gamma\delta- \beta(d-1)} \int_{A}
\int_{{\lambda}^{\beta} A} \int_0^{\infty} \int_0^{\infty} f(x,
h'_x{\lambda}^{-\gamma}) f({\lambda}^{-{\beta}} y' + x, h_y'
{\lambda}^{-\gamma} )\\
&&\phantom{{\lambda}^{2 - 2\gamma- 2\gamma\delta- \beta(d-1)}
\int_{A}
\int_{{\lambda}^{\beta} A} \int_0^{\infty} \int_0^{\infty}}
{}\times
c^{({\lambda})}_x ((\mathbf{0},h'_x), (y',h_y'))
\\
&&\phantom{{\lambda}^{2 - 2\gamma- 2\gamma\delta- \beta(d-1)}
\int_{A}
\int_{{\lambda}^{\beta} A} \int_0^{\infty} \int_0^{\infty}}
{}\times
{\rho}^{({\lambda})}(\mathbf{0},h'_x) {\rho}^{({\lambda
})}(y',h_y') \,dh'_x \,dh_y' \,dy' \,dx.
\end{eqnarray*}
Recalling from (\ref{twodef}) that
$\beta(d - 1) + \gamma(1 + {\delta}) = 1$, we have by definition of
$\tau$ [see (\ref{TAUU})] that $2 -
2\gamma- 2\gamma\delta- \beta(d-1) = 1 - \gamma(1+\delta) =
\tau$. After normalization by ${\lambda}^{\tau}$, the above integral
equals
%
\begin{eqnarray*}
&&\int_{A}
\int_{{\lambda}^{\beta} A} \int_0^{\infty} \int_0^{\infty} f(x,
h'_x{\lambda}^{-\gamma}) f({\lambda}^{-{\beta}} y' + x, h_y'
{\lambda}^{-\gamma} )\\
&&\phantom{\int_{A}
\int_{{\lambda}^{\beta} A} \int_0^{\infty} \int_0^{\infty}}
{}\times
g_{{\lambda}} (x,h'_x,y',h_y')\, dh'_x\,dh_y' \,dy' \,dx,
\end{eqnarray*}
where we put 
%
\[
g_{{\lambda}} (x,h'_x,y',h_y'):= c_x^{({\lambda})} ((\mathbf
{0},h_x'), (y',h_y'))
{\rho}^{({\lambda})}(\mathbf{0},h_x') {\rho}^{({\lambda})}(y',h_y').
\]

Clearly, $f(x, h'_x {\lambda}^{-\gamma}) f({\lambda}^{-{\beta}} y' +
x, h_y'
{\lambda}^{-\gamma} )$ converges to $f^2(x,0)$ as ${\lambda}\to
\infty$.
Lemma \ref{corlimit} implies for all $(x,h'_x,y',h_y') \in A_+
\times{\lambda}^{{\beta}} A \times\mathbb{R}_+$ that the product
$g_{{\lambda}}
(x,h'_x,y',h_y') (h'_x)^{-\delta} (h_y')^{-\delta}$ converges to
\[
c_x^{(\infty)}((\mathbf{0},h'_x),(y',h_y')) {\rho}^2_0(x)
\]
as ${\lambda}\to\infty$. Since, by Lemma \ref{corbds} and (R2),
$g_{{\lambda}} (x,h'_x,y',h_y') (h'_x)^{\delta} (h'_y)^{\delta}$ is
dominated in absolute value by the integrable function
\[
(x,h_x',y',h_y') \mapsto C' (h'_x)^{\delta} (h'_y)^{\delta} \exp
\biggl(-\frac{1}{C} \max\biggl(\frac{|y'| }{2},
h_x',h_y' \biggr) \biggr)
\]
on
$A_+ \times\mathbb{R}^{d-1} \times\mathbb{R}_+$, the dominated convergence
theorem combined with relation (\ref{SCALINGLIMIT2}) produces the
desired limit (\ref{varlimit}).

\subsection{\texorpdfstring{Proof of Theorem \textup{\protect\ref
{CLT}}}{Proof of Theorem 1.3}}\label{CLTsubsection}

Given Theorems 1.1 and 1.2, one may prove Theorem \ref{CLT} either
by the method of cumulants \cite{BY2} or by the Stein method
\cite{PY5}. The first approach shows that the Fourier transform
of ${\lambda}^{-\tau/2} \langle
f,\bar{\mu}^{\xi}_{{\lambda}\rho}\rangle$, namely,
\[
\mathbb{E}\exp[ i {\lambda}^{-\tau/2} \langle f,\bar
{\mu}^{\xi
}_{{\lambda}\rho}\rangle
],
\]
converges as ${\lambda}\to\infty$ to the Fourier transform
of a normal mean zero random variable with variance
$ \sigma_f^2
:= I(f^2) + J(f^2)$.
Even
though we use a formally different version of stabilization, this
is accomplished by following \cite{BY2} nearly verbatim. Indeed,
recall that Lemma \ref{corbds} shows the exponential decay of the
two point correlation function $c^{({\lambda})}_x((\mathbf{0},h'),
(y',h_y'))$.
In a similar way we may establish the exponential decay of
$k$-point correlation functions, and, more generally, that the
$k$-point correlation functions cluster exponentially, as shown in
Lemma 5.2 of \cite{BY2}. In this way we show (as in Lemma 5.3 of
\cite{BY2}) that for all $k = 3,4,\ldots$ and $f \in\mathcal
{C}_b(A_+)$ that
%
%
\begin{equation}\label{cumlimit}
\lim_{{\lambda}\to\infty} {\lambda}^{-\tau k/2} \langle
f^{\otimes k}, c_{\lambda}^k \rangle= 0,
\end{equation}
where $c_{\lambda}^k$ denotes the
$k$th cumulant figuring in the logarithm of the Laplace transform
(their existence follows by Lemma
\ref{expbds}). 
This
consequently shows that ${\lambda}^{-\tau/2} \langle
f,\bar{\mu}^{\xi}_{{\lambda}\rho}\rangle$ converges to a
mean zero
normal random variable with variance $\sigma_f^2$. The
convergence of the finite-dimensional distributions follows from
the Cram\'er--Wold device and is standard (see, e.g., page 251 of
\cite{BY2} or \cite{Pe1}).

Alternatively, we may also use the Stein method \cite{Pe1,PY5}.
This is a bit simpler and has the advantage of yielding rates of
convergence when $\sigma_f^2 > 0$, as would be the case when
$\delta= 0$ and $\alpha= 2$ (Lemma 7 of \cite{Re4} combined with
Section \ref{ApplSection} below) or when $\alpha= 1$ (Theorem 2.2
of \cite{BY4}). (When $\sigma_f^2 = 0$, then ${\lambda}^{\tau
/2}\langle
f, \bar{\mu}_{{\lambda}\rho}^{\xi} \rangle$ converges to a unit point
mass.) Our proof is based closely on \cite{PY5}, which uses a
formally different version of stabilization. For simplicity, we
assume $A = [0,1]^{d-1}$.


Recalling that $\bar{x} := (x, h_x)$, we have
\[
\langle f, \mu_{{\lambda}\rho}^{\xi} \rangle= \sum_{ \bar{x} \in
{\mathcal{P}}_{{\lambda}{\rho}}} \xi( \bar{x}, {\mathcal
{P}}_{{\lambda}{\rho}}) f (\bar{x}) =
\sum_{
\bar{x} \in{\mathcal{P}}_{{\lambda}{\rho}}} \xi^{({\lambda})}
\bigl(
(\mathbf{0}, h'_x), {\mathcal{P}}_{{\lambda}
{\rho}}^{({\lambda})}[x]\bigr) f ((x,h_x' {\lambda}^{-\gamma}) ).
\]
For all $L > 0$, let
\begin{eqnarray*}
T_{\lambda}&:=& T_{\lambda}(L):= \sum_{ \bar{x} \in{\mathcal
{P}}_{{\lambda}{\rho}} \cap
([0,1]^{d-1} \times[0,L {\lambda}^{-\gamma} \log{\lambda}]) } \xi
^{({\lambda})} \bigl(
(\mathbf{0}, h'_x), {\mathcal{P}}_{{\lambda}{\rho}}^{({\lambda
})}[x]\bigr) f ((x,h_x' {\lambda}^{-\gamma}) )
\\
&\hspace*{4pt}=& \sum_{ \bar{x} \in{\mathcal{P}}_{{\lambda}{\rho
}}: \ h_x' \leq
L {\lambda}^{-\gamma}
\log
{\lambda}} \xi^{({\lambda})} \bigl( (\mathbf{0}, h'_x), {\mathcal
{P}}_{{\lambda}{\rho}}^{({\lambda})}[x]\bigr) f
((x,h_x' {\lambda}^{-\gamma}) ).
\end{eqnarray*}

By Lemma 3.2, given arbitrarily large ${\kappa}> 0$, if $L$ is large
enough, then
$\langle f, \mu_{{\lambda}\rho}^{\xi} \rangle$ and $T_{\lambda}$ coincide
except on a set with probability $O({\lambda}^{-{\kappa}})$
in ${\lambda}$. Thus, $T_{\lambda}$ has the same
asymptotic distribution as $\langle f, \mu_{{\lambda}\rho}^{\xi}
\rangle$ and it suffices to find a rate of convergence to the
standard normal for $(T_{\lambda}- \mathbb{E}T_{\lambda})/\sqrt
{\operatorname{Var}T_{{\lambda}}}$.

Subdivide $[0,1]^{d-1}$ into $V({\lambda}):= {\lambda}^{{\beta}(d-1)}
({\rho}_{\lambda})^{-(d-1)}$ sub-cubes $C_i^{{\lambda}}$ of edge length
${\lambda}^{-\beta} {\rho}_{\lambda}$ and of volume ${\lambda
}^{-{\beta}(d-1)}
({\rho}_{\lambda})^{d-1}$, where $\rho_{\lambda}:= M \log{\lambda
}$ for some large
$M$, exactly as in Section 4 of \cite{PY5}.

Enumerate ${\mathcal{P}_{{\lambda}\rho} }\cap(C^{{\lambda}}_i
\times L
{\lambda}
^{-\gamma} \log{\lambda})$ by
$\{ \bar{X}_{i,j} \}_{j=1}^{N_i}$
where $\bar{X}_{i,j}:= (x_{ij},
h_{ij})$. Re-write $T_{\lambda}$ as
\[
T_{\lambda}= \sum_{i=1}^{V({\lambda})} \sum_{j = 1}^{N_i} \xi
^{({\lambda})} \bigl(
(\mathbf{0}, h'_{ij}), {\mathcal{P}}_{{\lambda}{\rho}}^{({\lambda
})}[x_{ij}]\bigr) f ((x_{ij},h'_{ij}
{\lambda}^{-\gamma} ) ).
\]
This is the analog of $T_{\lambda}$ in \cite{PY5}.

For any random variable $X$ and any $p > 0$, let
$\Vert X\Vert_p:= (\mathbb{E}[|X|^p ])^{1/p}.$
For all $1 \leq i \leq V({\lambda})$, we have $\sum_{j = 1}^{N_i}
\xi^{({\lambda})} ( (\mathbf{0}, h'_{ij}), {\mathcal{P}}_{{\lambda
}{\rho}}^{({\lambda})}[x_{ij}]) \leq
N_i$, where $N_i$ is Poisson with mean
\[
{\lambda}\int_{C^{{\lambda}}_i \times[0,L{\lambda}^{-\gamma} \log
{\lambda}]} \rho(u)
\,du = O([\log{\lambda}]^{1+\delta}).
\]
It follows by the boundedness of
$f$ that
%
\[
\Biggl\Vert\sum_{j = 1}^{N_i} \xi^{({\lambda})} \bigl( (\mathbf
{0}, h'_{ij}),
{\mathcal{P}}_{{\lambda}{\rho}}^{({\lambda})}[x_{ij}]\bigr)
f((x_{ij},h'_{ij} {\lambda}^{-\gamma} ))
\Biggr\Vert_3
\leq C \Vert f\Vert L^{1+\delta} (\log{\lambda})^{1+\delta} ({\rho
}_{\lambda})^{d-1},
\]
where $\Vert f\Vert$ denotes the essential supremum of $f$.
This is the analog of Lemma 4.3 in \cite{PY5} (putting $q =3$
there) with an extra logarithmic factor.

For all $1 \leq i \leq V({\lambda})$ and $j = 1,2,\ldots,$ let $R_{i,j}$
denote the radius of stabilization for $\xi^{({\lambda})}$ at
$X_{i,j}$ for ${\mathcal{P}}_{{\lambda}{\rho}}^{({\lambda})}$ if
$1 \leq j \leq N_i$ and
let $R_{i,j} $
be zero otherwise.

As in \cite{PY5}, put $E_{i}:= \bigcap_{j=1}^{\infty} \{R_{i,j}
\leq\rho_{\lambda}\}$ and let $E_{\lambda}:= \bigcap
_{i=1}^{V({\lambda})} E_i$.
Then by Lemma~\ref{StabLemma}(i), we have $P[E_{\lambda}^c] \leq
{\lambda
}^{-\kappa}$ for
$\kappa$ arbitrarily large if $M$ is large enough. This is the
analog of (4.11) of \cite{PY5}.

Next, recalling ${\rho}_{\lambda}= M \log{\lambda}$, we define the
analog of
$T'_{\lambda}$ in \cite{PY5}:
%
\[
T'_{\lambda}:= \sum_{i=1}^{V({\lambda})} \sum_{j = 1}^{N_i}
\xi^{({\lambda})}_{[{\rho}_{{\lambda}}]} \bigl( (\mathbf{0}, h'_{ij}),
{\mathcal{P}}_{{\lambda}
{\rho}}^{({\lambda})}[x_{ij}]\bigr) f ((x_{ij},h'_{ij} {\lambda
}^{-\gamma} ) ).
\]

Then we define, for all $1 \leq i \leq V({\lambda})$,
\[
S_i:= S_{Q_i}:= (\operatorname{Var}T'_{\lambda})^{-1/2} \sum_{j = 1}^{N_i}
\xi^{({\lambda})}_{[{\rho}_{{\lambda}}]} \bigl(
(\mathbf{0}, h'_{ij}), {\mathcal{P}}_{{\lambda}{\rho}}^{({\lambda
})}[x_{ij}]\bigr) \ f ((x_{ij},h'_{ij}
{\lambda}^{-\gamma} ) ).
\]
We define $S_{{\lambda}}:= \sum_{i=1}^{V({\lambda})} (S_i - \mathbb
{E}S_i)$, noting
that it is the analog of $S$ in \cite{PY5}.

Notice that $T'_{\lambda}$ is a close approximation to $T_{\lambda}$
and that,
by definition of $E_i, 1 \leq i \leq V({\lambda})$, it has a high amount
of independence between summands. 
In fact, by the independence property of Poisson point
processes, it follows that $S_i$ and $S_k$ are independent
whenever $d(C_i^{{\lambda}}, C_k^{{\lambda}})
> 2 {\lambda}^{-\beta} {\rho}_{{\lambda}}.$

Next we define a graph $G_{\lambda}:= (\mathcal{V}_{\lambda},
\mathcal{E}_{\lambda})$
as follows. The set $\mathcal{V}_{\lambda}$ consists of the sub-cubes
$C^{\lambda}_1,\ldots,C^{\lambda}_{V({\lambda})}$ and the edges
$(C^{\lambda}_i,C^{\lambda}_j)$
belong to $\mathcal{E}_{\lambda}$ if $d(C^{\lambda}_i, C^{\lambda}_j)
\leq2
{\lambda}^{-\beta} {\rho}_{\lambda}$. Since $S_i$ and $S_k$ are independent
whenever $d(C_i^{{\lambda}}, C_k^{{\lambda}}) > 2 {\lambda}^{-\beta
} {\rho}_{{\lambda}}$, it
follows that $G_{\lambda}$ is a dependency graph for
$\{S_i\}_{i=1}^{V({\lambda})}$.

Now proceed exactly as in \cite{PY5}, noting that:
\begin{longlist}
\item$V({\lambda}) = {\lambda}^{{\beta}(d-1)} ({\rho}_{\lambda
})^{-(d-1)}$,
\item the maximum degree of $G_{\lambda}$ is bounded by $5^d$,
\item for all $1 \leq i \leq V({\lambda})$, we have $\Vert S_i\Vert_3
\leq K
(\operatorname{Var}(T'_{\lambda}))^{-1/2} (\log{\lambda})^{1+\delta
}\times({\rho}_{\lambda})^{d-1} =:
\theta[{\lambda}]$.
\end{longlist}
%


These bounds correspond to the analogous bounds (i), (ii) and
(iii) on pages 54--55
of \cite{PY5}. Moreover, provided $\sigma_f^2 > 0$, then the
counterpart of (v) of \cite{PY5} holds, namely,
\[
\operatorname{Var}[T'_{\lambda}] = \Theta(\operatorname
{Var}[T_{\lambda}]) = \Theta({\lambda}^{\tau}).
\]

Putting $q = 3$ in (4.1) and
(4.18) of \cite{PY5} gives a rate of convergence for both
$S_{\lambda}$ and
$(T_{\lambda}- \mathbb{E}T_{\lambda})/\sqrt{\operatorname
{Var}T_{{\lambda}}}$ to the standard normal.
This rate is
\[
O( V({\lambda}) \theta[{\lambda}]^3) = O\bigl( {\lambda}^{{\beta}(d-1)}
({\rho}_{\lambda}
)^{-(d-1)} (
{\lambda}^{\tau})^{-3/2} (\log{\lambda})^{3(1+\delta)} \rho
_{\lambda}^{3(d-1)} \bigr).
\]
Recalling that $\tau= \beta(d-1)$, we rewrite this as
%
%
\begin{equation}
\label{rates}
O\bigl({\lambda}^{-\tau/2} \log{\lambda}^{3(1+\delta) +
2(d-1)} \bigr).
\end{equation}

This completes the proof of Theorem 1.3.\quad\qed



\section{Proofs of applications}\label{ApplSection}

The purpose of the present section is to derive Theorems 2.1 and 2.2
from our general theorems of Section \ref{GenRes}.

\begin{pf*}{Proof of Theorem \ref{convexhullthm}}
To derive Theorem \ref{convexhullthm}
from our general theory, we translate the convex hull problem into
the language of $\psi$-growth processes with overlap. To this end, recall
first that for a compact convex body $C \subseteq\mathbb{R}^d$ we
define its
support function $h_C \dvtx S_{d-1} \to\mathbb{R}$ by
\[
h_C(u) := \sup_{\bar{x} \in C} \langle\bar{x}, u \rangle
,\qquad u \in S_{d-1},
\]
with now $\langle\cdot, \cdot\rangle$ standing for the usual scalar
product in $\mathbb{R}^d$; see
Section 1.7 in \cite{SCHN}. An easily verified and yet crucial feature
of the support functional
$h_{ \cdot}(\cdot)$ is that
%
%
\begin{equation}\label{SuppFnct}
h_{\operatorname{conv}\{ \bar{x}_1,\ldots,\bar{x}_k \}}(u) = \max
_{1 \leq i \leq
k} h_{\{ \bar{x}_i \}}(u),\qquad u \in
S_{d-1},
\end{equation}
for each collection $\{ \bar{x}_1,\ldots, \bar{x}_k \}$ of points in
$\mathbb{R}^d.$
Moreover, by definition, it is clear that, for all $u \in S_{d-1}$,
we have $ h_{\{ \bar{x} \}}(u) = \langle\bar{x}, u \rangle$, $u \in
S_{d-1}$.

This leads to the following way of describing
$\mathcal{V}({\mathcal{P}_{{\lambda}\rho} })$ considered in
Theorem~\ref
{convexhullthm}.
For
a particular realization $\{ \bar{x}_1,\ldots, \bar{x}_k \}$ of
${\mathcal{P}_{{\lambda}\rho} }$
in $B_d$, we consider the collection $H[\bar{x}_1],\ldots,H[\bar
{x}_k]$ of
\textit{support epigraphs} given by
%
%
\begin{equation}\label{KDef}
H[\bar{x}] := \bigl\{ (y,h_y) \in S_{d-1} \times\mathbb{R}_+\dvtx
h_y \geq
1 -
h_{\{ \bar{x} \}}(y) \bigr\},
\end{equation}
where $h_y$ stands for the distance between $\bar y$ and the boundary
$S_{d-1} = \partial B_d.$
A~compact convex body is uniquely determined by its support functional
(cf. Section
1.7 in \cite{SCHN}),
and in view of (\ref{SuppFnct}), the set $\operatorname{conv}(\{ \bar
{x}_1,\ldots
,\bar{x}_k \})$
is in one-to-one correspondence with the union $\bigcup_{i=1}^k H[\bar{x}_i].$
Further, the number of vertices in the convex hull is easily seen to
coincide with the
number of those $\bar{x}_i$, $ i = 1,\ldots, k,$ for which $H[\bar
{x}_i]$ is
not completely contained in the union $\bigcup_{j \neq i} H[\bar{x}_j].$

Next we shall also write $r_y := 1 - h_y$ for the
distance between $\bar y$ and the origin of $\mathbb{R}^d$.
Note now that the
intensity measure
$\rho(\bar{x})\, d\bar{x}$, $\bar{x} \in B_d,$ coincides
with $\rho((x,r)) r^{d-1}\, dr\, dx = \rho((x,r)) (1-h)^{(d-1)}\,
dh\, dx$,
where $ \bar{x}:= (x,r),$
with $r \in[0,1]$ denoting the distance between $\bar{x}$ and the
origin of
$\mathbb{R}^d$,
with $h:=1-r$ and with $x \in S_{d-1}$ being the
radial projection of $\bar{x}$ onto $\partial B_d = S_{d-1}.$
Observe also that the support epigraph $ H[\bar{x}]$ as given in
(\ref{KDef}) can be represented by
\[
H[(x,r)] = \{ (y,h_y) \in S_{d-1} \times\mathbb{R}_+ \dvtx h_y \geq1
- r \cos
(\mathrm{dist}_{S_{d-1}}(x,y)) \}
\]
with
$\mathrm{dist}_{S_{d-1}}(x,y) := \cos^{-1} \langle x, y \rangle$ denoting
the geodesic distance in
$S_{d-1}$ between $x$ and $y.$ Now put
\[
\psi(l) := 1 - \cos(l).
\]
Writing the inequality $h_y \!\geq1 - r
\cos(\mathrm{dist}_{S_{d-1}}(x,y))$ as $h_y \geq1 - r +\break
r\psi(\mathrm{dist}_{S_{d-1}}(x, y))$, we have
%
%
\begin{equation}\label{KRepr}
H[(x,r)] = \{ (y,h_y) \in S_{d-1} \times\mathbb{R}_+ \dvtx h_y \geq h
+ r \psi
(\mathrm{dist}_{S_{d-1}}(x,y))
\},
\end{equation}
in other words, the support epigraphs are remarkably similar to
the upward cones (\ref{upcone}) described at the outset.

The above observations naturally suggest \textit{identifying the cardinality
of the studied set
$\mathcal{V}({\mathcal{P}_{{\lambda}\rho} })$ with the number of extreme
points in the $r\psi$-growth process with overlap in the sense of Section
\textup{\ref{GenRes}}
with the underlying point density} $\rho((x,r)) r^{d-1} =
\rho((x,r)) (1-h)^{d-1}.$
Likewise, the vertex empirical measure $\mu_{{\lambda}\rho}$ in
(\ref
{vepm}) corresponds
to the empirical measure $\mu^{\xi}_{{\lambda}\rho}$, $\xi:= \xi
(\psi;\cdot)$;
see (\ref{XiPsi}).

This
identification is valid modulo the following issues though:
\begin{longlist}
\item[(1)] the ``spatial'' coordinate $x$ of a point $\bar{x}:= (x,r)
\in B_d$
falls into $S_{d-1}$ rather than into a subset $A$ of $\mathbb{R}^{d-1},$
as required in Section \ref{GenRes},
\item[(2)] $\psi$ as given above is monotone only in a neighborhood
of $0,$
and moreover, we do not have $\lim_{l\to\infty} \psi(l) = \infty,$ which
violates $(\Psi1)$,
%
\item[(3)] the support epigraph $H[(x,r)]$ coincides with the
$(x,h)$-shifted\break\mbox{$\psi$-epigraph} $K[x,h]$ given by (\ref{upcone})
only when $r = 1$ and, hence, only when $h=0$; in general, for $0 \leq r
\leq1$, the set
$H[(x,r)]$ is an $(x,h)$-shifted $r \psi$-epigraph.
\end{longlist}

We claim, however, that the above three restrictions can be neglected
in the
asymptotic regime ${\lambda}\to\infty$, thus rendering the theory of
Section
\ref{GenRes}
applicable.
Indeed, first note that the sphere $S_{d-1}$, unlike the boundary
of a general smooth convex body, has a spatially homogeneous structure
and so the behavior of $\psi$ is independent of $x$, exactly as in
Section \ref{GenRes}. Moreover, the sphere $S_{d-1}$, being a smooth
manifold, has a local geometry
coinciding with that of $\mathbb{R}^{d-1}$, which takes care of issue
(1). Concerning issues (2) and (3), for each $r \in(0,1)$, the convex hull
$\operatorname{conv}({\mathcal{P}_{{\lambda}\rho} })$
coincides with $\operatorname{conv}({\mathcal{P}_{{\lambda}\rho}
}\cap
(B_d \setminus
B_d(0,r)))$ with
overwhelming probability, that is, the probability of the complement
event goes to zero exponentially fast in ${\lambda}$;
see the discussion in \cite{KUE} and the references therein.
This allows us to focus on the geometry of
$\operatorname{conv}({\mathcal{P}_{{\lambda}\rho} })$ in a thin shell
$B_d \setminus
B_d(0,r)$ within
a distance $1-r$ from the boundary $S_{d-1}.$

Consequently, only the behavior of $\psi$ in a
neighborhood of $0$ matters.
Recalling that the standard
re-scaling of Section \ref{ScaRe} involves scaling in the
spatial directions by ${\lambda}^{{\beta}}$, it follows that for a
given $\bar{x}:=
(x,r)$ and support epigraph $H[\bar{x}]$, the contribution of points
distant from $x$ by more than $O({\lambda}^{-\beta})$ is negligible
in view of the argument in Lemma \ref{StabLemma}(i) and no distortions
from the local Euclidean geometry have to be taken into account
in the limit under this re-scaling.
Likewise, we only
have to control the geometry
of $H[\bar{x}]$, $\bar{x} := (x,r),$ for $r$ arbitrarily close to~$1$.
This allows us to rewrite the proofs of Theorems \ref{LLN}--\ref{CLT}
for the thus modified $r$-dependent $\psi.$ Indeed, the stabilization Lemma
\ref{StabLemma}, as well as Lemma \ref{expbds}, do not require any
modifications in their proofs and neither does Lemma 3.4 nor Lemma 3.5.
Consequently, the arguments leading to the central limit theorem in Section
\ref{CLTsubsection} do not require modification either. In this
context we note that the proof of Lemma \ref{StabLemma} would
break down if the sphere $S_{d-1}$ were replaced by
a nonconvex set allowing for long-range dependencies between
extreme points.

It only remains to show the limit arguments in Sections \ref{VARproof}
and \ref{LLNproof}
remain valid for the modified $\psi.$ To see that this is indeed the
case, we note that the arguments rely on two main ingredients: on
stabilization which holds with no changes as stated above, and on
re-scaling relations discussed in Section \ref{ScaRe}. However,
it is easily seen that the re-scaling relations and their proofs
can be readily rewritten for the modified $\psi$, the only essential
modification being to add one extra argument ($h := 1-r$) to the $\psi
$-function,
which anyway vanishes in the scaling limit of Section \ref{ScaRe}
with $h = 1-r$ tending to $0$ as discussed above (whereas the contribution
coming from smaller $h$ is negligible in view of Lemma \ref{expbds}).
This discussion takes care of issues (2) and (3) above.

Thus, we can now conclude that the considered convex hull process falls
into the range of
applicability of the general theory of Section
\ref{GenRes}, with $\alpha= 2$ in
$(\Psi2)$ and $\delta$ in (R2) coinciding with that in
the statement of
Theorem \ref{convexhullthm}. Thus, we obtain the required Theorem \ref
{convexhullthm}
as a consequence
of the general Theorems \ref{LLN}--\ref{CLT}.
The rate of
convergence follows from $(\ref{rates})$ by putting $\delta= 0$
and $\alpha= 2$.
\end{pf*}
%

\begin{pf*}{Proofs of Theorems \ref{maximalthm} and \ref{maximalthmBin}}
Theorem \ref{maximalthm} follows directly by the general theory in
Section \ref{GenRes} (Theorems \ref{LLN}--\ref{CLT} with $\alpha\in
(0,1]$).
The rate (\ref{maxptrate}) follows from (\ref{rates}) by putting
$\delta= 0$ and $\alpha= 1$. We thus focus attention on
establishing Theorem \ref{maximalthmBin}.
The first lemma yields (\ref{explimmaxBin}).
\noqed
\end{pf*}


%
\begin{lemm}\label{lem4.1} 
For all $f \in\mathcal{C}_b(A_+)$, we have
%
%
\begin{equation}
\label{differ}
|\mathbb{E}[\langle f, {\nu}_n^{\xi} \rangle] -
\mathbb{E}[\langle
f, \mu_{n \rho}^{\xi} \rangle] | = O(n^{-{\tau}'}).
\end{equation}
\end{lemm}


\begin{pf}
For all $\bar w \in A_+$, let $p(\bar w):=
\int_{K^{\downarrow} [\bar w]} {\rho}(u) \,du,$
where
$ K^{\downarrow}[\bar{w}]$ is as in (\ref{downarrow}) with
$\psi(l) = l^{\alpha}.$
Note that in our current setting for all $w
\in A_+$ we have $p(w) \in[0,1]$ since
${\rho}$ is a probability density. Also, note that
$\psi^{(n)}
\equiv\psi$ and $K^{(n)} \equiv K$ with $K^{(n)} :=
\{ (y^{(n)},h_y^{(n)})\dvtx(y,h_y) \in K \},$
that is, the self-similarity
under the re-scaling is immediate rather than emerging as $n \to\infty.$
For all $s
\in[0,1]$
and $f \in\mathcal{C}_b(A_+)$, let $B_f(s):= \int_{p(\bar w) \leq s}
f(\bar w)
{\rho}(\bar w)\, d \bar w$.
Recalling that for $\alpha\in(0,1]$ the \mbox{$\psi$-extremality}
of a
point $w$ in a given sample is equivalent to having no other
sample points in $K^{\downarrow}[w]$ (see the discussion at the
beginning of Section \ref{NMAX}), we have
\begin{eqnarray*}
\mathbb{E}[\langle f,{\nu}_n^{\xi} \rangle] &=& n \int_{A_+} \bigl
(1 - p(w)\bigr)^{n-1}
f(w) {\rho}(w)\, dw
\\
&=& n \int_0^1 \int_{p(w) = s} (1 -
s)^{n-1}f(w)\rho(w) \,dw \,ds = n \int_0^1 (1 - s)^{n - 1}\, dB_f(s)
\end{eqnarray*}
by Fubini's theorem. Similarly,
%
%
\begin{equation}
\label{tauber} \mathbb{E}[\langle f,
\mu_{n \rho}^{\xi} \rangle] = n \int_0^1 e^{-ns} dB_f(s) \sim C_f
n^{\tau},
\end{equation}
where the asymptotics are given by Theorem \ref{LLN}. $B_f$
is monotone, nondecreasing and Karamata's Tauberian theorem (e.g.,
Theorem 2.3 in \cite{Se}) gives $B_f(s) \sim C_f s^{{\tau}'}$ as
$s \to0^+$.
Notice
\begin{eqnarray*}
| \mathbb{E}[\langle f, \mu_{n \rho}^{\xi} \rangle] -\mathbb
{E}[\langle f,
{\nu}_n^{\xi} \rangle] | &=& n \int_0^1 \bigl(e^{-ns} - (1 - s)^{n
-1}\bigr)\, dB_f(s)
\\
&\leq& n \int_0^1 \bigl(e^{-ns} - e^{n \ln(1 - s)}\bigr)\, dB_f(s)
\\
&\leq& Cn^2 \int_0^1 e^{-ns} s^2 \,dB_f(s)
\\
&=& Cn^2 \int_0^{1/n} e^{-ns} s^2\, dB_f(s) + Cn^2 \int_{1/n}^1
e^{-ns} s^2 \,dB_f(s).
\end{eqnarray*}
The first integral behaves
like $Cn^{-{\tau}'}$ since $B_f(s) \sim C_f s^{{\tau}'}$, whereas
the second behaves like ${C\over n} \int_1^n u^2 e^{-u} dB_f(u/n)
\leq C/n$, since $B_f$ is bounded by $B_f(1).$
This gives (\ref{differ}).
\end{pf}

We now establish the remainder of Theorem
\ref{maximalthmBin}.
Recall $\bar{u}' := (u',h'_u)$.
For all ${\lambda}> 0$, define
\[
A'({\lambda}) := \biggl\{ \bar{y}' \in{\lambda}^\beta A \times
\mathbb{R}_+ \dvtx\int_{
K^{\downarrow} [\bar{y}']} \rho^{({\lambda})} (u) \,du \leq C \log
{\lambda} \biggr\}.
\]
Let $A({\lambda}):= \{ \bar{y} \in A_+: \bar{y}' \in A'({\lambda})
\}$ and
put $a_{\lambda}:=\int_{A({\lambda})} \rho(w) \,dw$. Note that by Lemma
\ref{expbds} the probability that a sample point from ${\bar{\mathcal{X}_n}
} := \{ X_i \}_{i=1}^n$ in $A_+ \setminus A({\lambda})$ is
$\psi$-extremal is at most
%
%
\begin{equation}\label{P25}
C \exp\biggl(- \frac{ n \log{\lambda}}{C {\lambda}} \biggr)
\end{equation}
and the same holds for $\bar{\mathcal{X}}_n$ replaced by the Poisson sample
with intensity $n\rho.$ Indeed, although Lemma~\ref{expbds} was
originally established for Poisson samples, it is easily seen that
the same proof works also for binomial samples, as it essentially
relies on exponentially decaying upper bounds for probabilities of
certain sets in $A_+$ being devoid of points of the underlying
point process. 
Thus, the $\psi$-extremal points are predominantly concentrated in
$A({\lambda})$, a fact which we will use to show (\ref{varlimmaxBin}).
First we find growth bounds for $a_{\lambda}$.


\begin{lemm}\label{LMA} We have
$a_{\lambda}\leq C ( \log{\lambda})^{ \alpha( 1 + {\delta})/
(\alpha+ d - 1)
} {\lambda}^{ -\alpha(1 + {\delta})/(d - 1 + \alpha(1 +
\delta)) }$.
\end{lemm}

\begin{pf}
If $M({\lambda}):= \sup\{h_y\dvtx h_y \in A({\lambda})\}$,
then note
that $\alpha({\lambda})$ grows like $\int_0^{M({\lambda})}
h_y^{\delta}\, dh_y =
C(M({\lambda}))^{1 + {\delta}}$. We now find $M({\lambda})$.

If $\bar{y}':=(y',h_y') \in A'({\lambda})$, then, by (\ref
{volinvert}), we
have ${h'}_y^{(\alpha+ d - 1)/\alpha} \leq C \log{\lambda}.$
Since $h'_y = {\lambda}^{\gamma} h_y$ and since $\gamma= \beta
\alpha$,
it follows that
\[
h_y^{(\alpha+ d - 1)/\alpha} {\lambda}^{ \beta(\alpha+ d - 1)
} \leq C \log{\lambda}.
\]
Since $\gamma(d-1)/\alpha= \tau$, we have
\[
h_y^{(\alpha+ d - 1)/
\alpha} {\lambda}^{ \gamma+ \tau} \leq C \log{\lambda}\quad\mbox{or}
\quad h_y
\leq( \log{\lambda})^{\alpha/(\alpha+ d - 1) } {\lambda}^{
-\alpha( \gamma+ \tau)/(\alpha+ d - 1) },
\]
that is,
\[
M({\lambda}) \leq( \log{\lambda})^{\alpha/(\alpha+ d - 1)}
{\lambda}^{ {-\alpha( \gamma+ \tau)}/(\alpha+ d - 1) } .
\]

Recall that $\gamma+ \tau= (\alpha+ d - 1)/(d-1 + \alpha(1 +
\delta))$ to get the result. 
\end{pf}

The next lemma yields (\ref{varlimmaxBin}). The proof borrows
heavily from \cite{BY4} and, for the sake of completeness, we
provide the details.

\begin{lemm}\label{varlimB}
For all $f \in\mathcal{C}_b(A_+)$, we have
$ \lim_{n \to\infty} n^{-\tau} \operatorname{Var}[\langle f, {\nu
}_n^{\xi}
\rangle] = \lim_{n \to\infty} n^{-\tau} \operatorname
{Var}[\langle f, \mu_{n
\rho}^{\xi} \rangle].$
\end{lemm}
\begin{pf}
Recall $\bar{\mathcal{X}_n}:= \{X_i\}_{i=1}^n$.
Let $N_n:=\mbox{card} \{
{\bar{\mathcal{X}_n} } \cap A(n) \}$ and $N'_n:=\mbox{card} \{
{\mathcal{P}}_{n{\rho}}
\cap A(n) \}$. For all $r = 1,2,\ldots,$ denote by $e(r):= e_f(r)$
the expected value of the functional $\langle f\cdot
\bold{1}(A(n)), \mu_n^{\xi} \rangle$ \textit{conditioned on
$\{N(n)=r\}$}, and by $v(r)$ the variance of this functional
conditioned on $\{N(n) = r\}$.
Let ${\nu}_n^A := {\nu}_n^{\xi,A}$ denote the point
measure induced by the $\psi$-extremal points in $\{ {\bar{\mathcal
{X}_n} }\cap A(n) \}$. Similarly, let
$\mu_{n\rho}^A := \mu_{n \rho}^{\xi,A}$ denote the point
measure induced by the $\psi$-extremal points in $ \{ {\mathcal
{P}}_{n{\rho}} \cap A(n) \}$.

By the bound (\ref{P25}) on the probability of a given point
outside $A(n)$ being extremal,
${\nu}_n^A$ coincides with ${\nu}_n$ and $\mu_{n\rho}^{\xi,A}$
coincides with $\mu^{\xi}_{n\rho}$ except on a set with
probability at most $n C \exp( - C \log n) $ $= C n^{-C+1}.$ Since
$C$ can be chosen arbitrarily large, it suffices to show that
%
%
\begin{equation}
\label{equalvar} \lim_{n \to\infty} n^{-\tau}\operatorname
{Var}[\langle f,
{\nu}_n^{A} \rangle] = \lim_{n \to\infty} n^{-\tau}\operatorname
{Var}[\langle f,
\mu_{n \rho}^{A} \rangle].
\end{equation}

The conditional variance formula implies that
\begin{eqnarray*}
\operatorname{Var}[\langle f, {\nu}_n^A\rangle] &=&\operatorname
{Var}[e(N_n)] + \mathbb{E}[ v(N_n)] \quad
\mbox{and} \\
\operatorname{Var}[\langle f,\mu_{n \rho}^{A}
\rangle] &= &\operatorname{Var}
[e(N'_n)] + \mathbb{E}[v(N'_n)].
\end{eqnarray*}

We prove (\ref{equalvar}) by showing that:
\begin{longlist}
\item the terms $\mathbb{E}[ v(N_n)] $ and $\mathbb{E}[v(N'_n)]$ are
dominant and
that their ratio tends to one as $n \to\infty$, and
\item$\operatorname{Var}[e(N_n)]$ and $\operatorname{Var}[e(N'_n)]$
are both $o(n^{\tau})$.
\end{longlist}

We will first show (ii) as follows. For all $s > 0$, recall that
$B_f(s):= \int_{p(\bar w) \leq s} f(\bar w) {\rho}(\bar w)\, d \bar w.$
By Fubini's theorem, for all $r = 1,2,\ldots$ and with $a_n
= \int_{A(n)} {\rho}(w) \,dw$, 
we obtain
\[
e(r)= \frac{r}{a_n} \int_{A(n)} \biggl(1- \frac{p(w)}{a_n}
\biggr)^{r-1} f(w) {\rho}(w) \, dw = \frac{r}{a_n } \int_{0}^{a_n}
\biggl(1- \frac{s}{a_n} \biggr)^{r-1} \,dB_f(s).
\]
Letting $\Delta_r$ denote the difference $e(r+1)-e(r)$, we obtain
\[
\Delta_r = \frac{1}{a_n} \int_0^{a_n} \biggl(1-\frac{s}{
a_n} \biggr)^{r} - \frac{rs}{ a_n} \biggl(1-\frac{s}{
a_n} \biggr)^{r-1}\,dB_f(s).
\]
Setting $u = rs/ a_n$ and
applying $B_f(s) \sim C_f s^{\tau'}$, we see that ($\tau= 1 -
{\tau}'$)
\[
|\Delta_r| \leq\frac{C_f }{r} \int_0^{r} \biggl| \biggl(1-\frac{u}{
r} \biggr)^r - u \biggl(1-\frac{u}{r} \biggr)^{r-1} \biggr|
\biggl(\frac{ua_n }{r} \biggr)^{ -\tau} \,du.
\]
Since $ \sup_{r > 0} \int_0^{r} | (1-\frac{u}{r} )^r
- u (1-\frac{u}{r} )^{r-1} | u^{-\tau} \,du \leq C$,
it follows that $|\Delta_r| \leq C
(\frac{a_n }{r} )^{ -\tau}$ .

When $r \in I_n:= (na_n - C (\log n) (na_n)^{1/2}, na_n
+ C (\log n) (na_n)^{1/2})$, then, by Lemma \ref{LMA}, for $n$
large,
\begin{eqnarray*}
|\Delta_r| &\leq& C (na_n)^{-1} n^{ \tau} = C a_n^{-1} n^{
-{\tau}'} \\
&=& C n^{ \alpha(1 + {\delta})/(d-1 + \alpha(1 +
\delta))
} n^{ -{\tau}'} (\log n)^{ -\alpha(1 + {\delta})/(\alpha+ d -
1) }.
\end{eqnarray*}
Recalling that ${\tau}' = (1 + {\delta})\alpha/(d-1 + \alpha(1 +
\delta))$, we see that for $r \in I_n$ we have
\[
|\Delta_r| \leq\Delta(n) := C (\log n)^{ -\alpha(1 + \delta)
/(\alpha+ d - 1) }.
\]

Write $e(N_n) = e(1) + \sum_{j = 2}^{N_n}(e(j) - e(j - 1))$ and
observe that $e(N_n)$ differs from the constant
$e(1) + \sum_{j = 2}^{E[N_n] }(e(j) - e(j - 1))$ by at most
\[
\sum_{j \in J_n} \bigl(e(j) - e(j - 1)\bigr),
\]
where $J_n := ( \min( \mathbb{E}[N_n], N_n ), \ \max( \mathbb
{E}[N_n], N_n )
)$. Thus,
\begin{eqnarray*}
\operatorname{Var}[e(N_n)] &\leq&\mathbb{E} \Biggl[ \sum_{j \in J_n}
\bigl(e(j) - e(j - 1)\bigr)
\Biggr]^2\\
& \leq&\mathbb{E} \Biggl[ \sum_{j \in J_n} \bigl(e(j) - e(j -
1)\bigr)
\mathbf{1}_{N_n \in I_n} \Biggr]^2 + o(1),
\end{eqnarray*}
by Cauchy--Schwarz and since (by increasing $C$ in the definition
of $I_n$) standard concentration inequalities (see, e.g.,
Proposition A.2.3(ii), (iii) and Proposition A.2.5(ii), (iii) in
\cite{BHJ}) show that $P[N_n \in I_n^c]$ can be made smaller than
any negative power of $n$.

For $j \in J_n$ and
$N_n \in I_n$, we have $j \in I_n$ and so $(e(j) - e(j - 1)) \leq
\Delta(n)$. Since the length of $J_n$ is bounded
by $|N_n - \mathbb{E}N_n|$, it follows that $\operatorname
{Var}[e(N_n)] \leq\operatorname{Var}[N_n]
(\Delta(n))^2 + o(1)$. Note that $\operatorname{Var}[N_n] \leq
Cn^{\tau}(\log
n)^{\alpha(1 + {\delta})/ (\alpha+ d - 1)}$. It follows that
$\operatorname{Var}[e(N_n)] \leq Cn^{\tau}(\log n)^{-\alpha(1 +
{\delta})/ (\alpha+
d - 1)} + o(1)$, that is, $\operatorname{Var}[e(N_n)] = o(n^{\tau}).$
Similarly,
$\operatorname{Var}[e(N'_n)] = o(n^{\tau})$ and so condition (ii) holds.

We now show condition (i) by showing that the ratio $\mathbb{E}
[v(N_n)]/\mathbb{E}[v(N'_n)]$ is asymptotically one, as $n\to\infty
$. Let
$p_{n,r}:=P[N_n=r]$ and $p'_{n,r}:=P[N'_n=r]$. Stirling's formula
implies that, for $|r-a_n n|\leq n^\beta$, where $0 < \beta< 1/2$,
%
%
\begin{equation}\label{second} \lim_{n \to\infty} \frac
{p_{n,r}}{p'_{n,r}} =
1
\end{equation}
uniformly. Now, for $|r-a_n n|> n^\beta$, where $\beta\in
(0,1/2)$
is chosen so that $n^{2 \beta}/na_n$ grows faster than some (small)
power of $n$,
we have that both $p_{n,r}$ and
$p'_{n,r}$ are bounded by $C\exp(-n^{{\delta}}/C)$ for some $C, \
{\delta}>
0$ (see, e.g., Proposition A.2.3(i) and Proposition A.2.5(i) in
\cite{BHJ}).
Write
\[
\mathbb{E}[v(N_n)]=\sum_{|r- a_n n| \leq n^\beta} v(r) p_{n,r}+ \sum_{|r-
a_n n|> n^\beta} v(r) p_{n,r}.
\]
The second sum is negligible since $0 < v(r) < r^2$ and $p_{n,r}$
is exponentially small. Consider the terms in the first sum. By
(\ref{second}), we have
$p_{n,r} = p'_{n,r}(1 + o(1))$ uniformly for all $|r - a_n| \leq
n^\beta$
and since the terms in the first sum are positive, it follows that
%
%
\begin{equation}\label{third} \lim_{n \to\infty}
\frac{\mathbb{E}[v(N_n)]} {\mathbb{E}[ v(N_n')] } = 1.
\end{equation}

Now from before we know $\operatorname{Var}[ \langle f, \mu_{n \rho
}^A \rangle]$
has asymptotic
growth $Cn^{\tau}, C > 0$. It follows that $\mathbb{E}[v(N_n')]$ has the
same growth, since $\operatorname{Var}[e(N_n')] =
o( n^{\tau} ).$ Thus, by (\ref{third}) and the
growth bounds $\operatorname{Var}[ e(N_n)] = o( n^{\tau} )$ and
$\operatorname{Var}[e(N_n')]
= o(n^{\tau})$, the desired identity (\ref{equalvar}) follows,
completing the proof of Lemma \ref{varlimB}.
\end{pf}

We conclude the proof of Theorem \ref{maximalthmBin} by
showing for all $f \in\mathcal{C}_b(A_+)$
\[
\lim_{n \to\infty} d_{\mathrm{TV}}( n^{-\tau/2} \langle f,
\bar{{\nu}}_n^{\xi} \rangle,\
n^{-\tau/2}\langle f, \bar{\mu}_{n \rho}^{\xi} \rangle) = 0,
\]
where the total variation distance between two measures $m_1$ and
$m_2$ is\break $d_{\mathrm{TV}}(m_1, m_2):= \sup_B |m_1(B) -
m_2(B)|$, where the
sup runs over all Borel subsets in $\mathbb{R}^d$. Since $
n^{-\tau/2} | \mathbb{E}[\langle f, {{\nu}}_n^{\xi}\rangle] -
\mathbb{E}[\langle f,
{\mu}_{n \rho}^{\xi} \rangle] | \to0$ by (\ref{differ})
and since $
n^{-\tau/2}\langle f, \bar{\mu}_{n \rho}^{\xi} \rangle$
converges in law to an appropriate Gaussian distribution,
recalling that $a_n = o(1)$ (see Lemma \ref{LMA}), Theorem
\ref{maximalthmBin} follows at once from the following:

\begin{lemm}
For all $f \in\mathcal{C}_b(A_+)$, we have
%
%
\begin{equation}\label{TV}
d_{\mathrm{TV}}( \langle f,{{\nu}}_n^{\xi} \rangle, \langle f, {\mu
}_{n \rho
}^{\xi}\rangle) = O(a_n).
\end{equation}
\end{lemm}

\begin{pf}
We follow the proof of Lemma 7.1 in \cite{BY4}.
Recall that ${\nu}_n^A$ is the measure induced by the maximal points
in $\{(X_i,h_i)\}_{i=1}^n \cap A(n)$ and, similarly, let $\mu_{n
\rho}^{\xi,A}$ be the measure induced by the maximal points in
${\mathcal{P}}_{n {\rho}} \cap A(n)$. If $C$ is large enough in the
definition of
$A(n)$, then the probability that points in $A_+ \setminus A(n)$
contribute to ${\nu}_n^{\xi}$ or $\mu_{n \rho}^{\xi}$ is
$O(n^{-2})$. It follows that, for all $f \in\mathcal{C}_b(A_+)$
\[
d_{\mathrm{TV}}( \langle f, \mu_{n \rho}^{\xi} \rangle, \langle f,
\mu_{n
\rho}^{\xi,A} \rangle) = O(n^{-2}) = o(a_n)
\]
and
\[
d_{\mathrm{TV}}( \langle f, {\nu}_n^{\xi} \rangle, \langle f, \mu
_n^{\xi,A}
\rangle) = O(n^{-2}) = o(a_n).
\]
Thus, we only need to show $d_{\mathrm{TV}}( \langle f, {\nu}_n^A
\rangle,
\langle f, {\nu}_{n \rho}^A \rangle) = O(a_n). $

Recall that $N_n$ is the number of points from ${\bar{\mathcal{X}_n} }$
belonging to $A(n)$.
Conditional on $N = r$, $\langle f, {\nu}_n^A \rangle$ is
distributed as $\langle f, \tilde{{\nu}}_r^A \rangle$, where
$\tilde{{\nu}}_r^A$ is the point measure induced by considering the
maximal points among $r$ points placed randomly according to the
restriction of $\rho$ to $A(n)$. The same is true for $\langle f,
\mu_{n \rho}^{\xi,A} \rangle$ conditional on the cardinality of
$\{{\mathcal{P}}_{n \rho} \cap A(n) \}$ taking the value $r$.

Hence, with $\mathit{Bi}(n,p)$ standing for a binomial random variable with
parameters $n$ and $p$ and $\mathit{Po}(\alpha)$ standing for a Poisson
random variable with parameter $\alpha$, we have for all $f \in
\mathcal{C}_b(A_+)$
\[
d_{\mathrm{TV}}( \langle f, {\nu}_n^A \rangle, \langle f, \mu_{n
\rho}^A
\rangle) \leq C d_{\mathrm{TV}} (\mathit{Bi} (n, a_n ), \ \mathit
{Po}(n a_n) )
\leq C \frac{1} {n a_n } \sum_{i=1}^n ( a_n)^2 \leq C a_n,
\]
where the penultimate inequality follows by standard Poisson
approximation bounds (see, e.g., (1.23) of Barbour, Holst and
Janson \cite{BHJ}). This is the desired estimate (\ref{TV}).
\end{pf}

\section*{Acknowledgments}
The authors gratefully
acknowledge helpful discussions with Yuliy Baryshnikov, who, in
particular, encouraged developing theory under the general
condition (R2) and who also contributed to the proofs of
Lemmas \ref{lem4.1} and \ref{varlimB}. The authors also thank an
anonymous referee
for comments leading to an improved exposition.

\printaddresses

\end{document}